\documentclass[12pt]{article}
\usepackage{amsfonts,amsmath,latexsym,amssymb,mathrsfs, url}
\usepackage{graphicx}
\usepackage{xcolor}
\usepackage{subcaption}

\numberwithin{equation}{section}

\def\numberlikeadb{\global\def\theequation{\thesection.\arabic{equation}}}
\numberlikeadb
\newcommand{\Cov}{{\mbox{Cov}}}
\newcommand{\beas}{\begin{eqnarray*}}
\newcommand{\enas}{\end{eqnarray*}}
\newcommand{\eqs}{\begin{eqnarray*}}
\newcommand{\ens}{\end{eqnarray*}}

\newcommand{\eqa}{\begin{eqnarray}}
\newcommand{\ena}{\end{eqnarray}}
\newcommand{\eq}{\begin{equation}}
\newcommand{\en}{\end{equation}}

\newcommand{\by}{\mathbf{y}}

\def\ignore#1{}

\def\half{{\textstyle{\frac12}}}

\def\quarter{{\textstyle{\frac14}}}

\def\urho{^{(\r)}}
\def\ul{^{(l)}}
\def\uld{^{(l')}}

\def\Ref#1{(\ref{#1})}

\def\b{\beta}
\def\s{\sigma}

\def\l{\lambda}

\def\sjn{\sum_{j=1}^n}

\def\hmu{{\hat\mu}}
\def\hnu{{\hat\nu}}
\def\hD{{\widehat\Delta}}
\def\hQ{{\widehat Q}}
\def\htau{{\hat\tau}}

\def\Def{\ :=\ }

\def\re{\mathbb{R}}

\def\tfrac#1#2{{\textstyle{\frac#1#2}}}

\def\non{\nonumber}
\def\th{\theta}

\def\e{\varepsilon}
\def\m{\mu}

\def\var{{\rm Var}}

\def\r{\rho}

\def\t{\tau}
\def\g{\gamma}

\def\ui{^{(1)}}
\def\um{^{(m)}}
\def\ut{^{(2)}}
\def\tu{{\tilde u}}

\def\nin{\noindent}
\def\D{\Delta}

\def\ex{{\mathbb E}}

\def\msk{\medskip}
\def\bsk{\bigskip}

\def\Eq{\ =\ }
\def\Le{\ \le\ }

\def\adbp{}

\def\tH{{\widetilde H}}
\def\G{\Gamma}

\def\ignore#1{}

\def\ssc{\psi}
\def\ssch{{\hat \ssc}}

\def\cupdot{\cup\kern-8.2pt\cdot\kern5.5pt}
\def\cov{\Cov}

\def\tX{\widetilde X}

\def\bone{{\bf{1}}}

\def\Def{\ :=\ }

\def\htau{{\hat{\tau}}}

\def\p{\pi}

\def\beps{\boldsymbol{\e}}
\def\bnu{\boldsymbol{\nu}}
\def\MVN{{\mathrm{MVN}\,}}

\def\TT{^\top}
\def\rhat{{\hat\r}}
\def\shat{{\hat\s}}
\def\Tr{{\mathrm{Tr}}}
\def\sn{\sum_{i=1}^n}
\def\cN{{\mathcal N}}
\def\ts{{\tilde s}}
\def\tbe{{\tilde \b}}
\def\parg{\g}
\def\hrho{\rhat}
\def\bard{{\overline d}}
\def\bhat{{\hat\b}}

\def\hbe{{\hat{\e}}}
\def\tbnu{{\bf \tilde\bnu}}
\def\adbb{}
\def\adb{}
\def\red{}

\def\EEE{L}
\def\ML{_{\mathrm{ML}}}
\def\VV{S}
\def\lmin{\l_{\mathrm{min}}}
\def\lmax{\l_{\mathrm{max}}}

\begin{document}

\title{Estimating the correlation in network disturbance models}

\author{
A. D. Barbour\footnote{Institut f\"ur Mathematik, Universit\"at Z\"urich,
Winterthurertrasse 190, CH-8057 Z\"URICH.
\msk}
\ and
Gesine Reinert\footnote{Department of Statistics,
University of Oxford, 24--29 St Giles', OXFORD OX1 3LB, UK.
GDR was supported in part by EPSRC grant EP/R018472/1 and by the COSTNET COST Action CA 15109.
}\\
Universit\"at Z\"urich and University of Oxford }

\date{}
\maketitle

\begin{abstract}
The Network Disturbance Model of \cite{doreian1989network} expresses the dependency between
observations taken at the vertices of a network by modelling the correlation
\adbb{between neighbouring vertices}, using
a single correlation parameter~$\rho$.  It has been observed that estimation of~$\rho$ in
{\it dense\/} graphs, using the method of Maximum Likelihood, leads to results that can
be both biased and very unstable.  In this paper, we sketch why this is the case, showing that
the variability cannot be avoided, no matter how large the network.  We also propose a more
intuitive estimator of~$\rho$, which shows little bias.
\adbb{The related Network Effects Model is briefly discussed.}
\end{abstract}

\nin {\bf Keywords:}  Network disturbance, network autocorrelation, maximum likelihood.

\msk\nin
{\bf MRC subject classification:} 91D30; 91Cxx, 62P25, 62J05

\section{Introduction}\label{intro}
 \setcounter{equation}{0}
The classical Gaussian linear regression model can be written as
\begin{equation}\label{linear-regression}  
    y \Eq  X \b + \e.
\end{equation} 
Here, the $n$-dimensional vector~$y$ of data is described as the sum of a structural element
$X \b$, where~$X$ is a known $n\times m$ matrix of full rank~$m$, 
$\b = (\b_1,\ldots,\b_m)\TT$ is an $m$-vector of unknown parameters, to be estimated,
and the elements $\e_1,\ldots,\e_n$ of the $n$-vector~$\e$ are realizations of independent
normally distributed random variables~$\beps_1,\ldots,\beps_n$, 
each with zero mean and variance~$\s^2$.  Here and throughout,
we conventionally use bold face to denote random elements, and italic for their realizations.

If the indices $i = 1,2,\ldots,n$ represent the vertices of \adbb{a \adbb{simple} undirected graph, whose 
(symmetric) adjacency} matrix is denoted by~$A$, it may be reasonable to suppose
that the errors~$\beps_i$ and~$\beps_j$ at neighbouring vertices $i$ and~$j$ are correlated, instead of being independent. 
A similar idea lies at the basis of time series analysis, where the error at a given time point
is assumed to depend on the errors at the closest previous time points.  In analogy with the
moving average models of time series analysis, the {\it Network Disturbance Model\/} 
of \cite{doreian1989network} (see also \cite[(4.6)]{ord1975}) supposes that
$\by \Eq  X \b + \beps$, where the errors~$\beps$ are now implicitly given by the equation
\eq\label{NDM-errors}
          \beps_i \Eq \rho\sum_{{j=1}}^n W_{ij}\beps_j + \bnu_i,\qquad 1\le i\le n,
\en
\adbb{or more succinctly, by
\eq\label{NDM-errors-2}
    K_\r \beps \Eq \bnu, \qquad\mbox{where}\quad K_\r \Def I - \r W.
\en
}
Here, the weight matrix~$W$ is known, the parameter~$\rho$ is unknown, and the random variables $\bnu_1,\ldots,\bnu_n$ 
are independent normal random variables with mean zero and variance~$\s^2$.  
The weight~$W_{ij} \ge 0$ can be thought of 
as a (relative) measure of
the influence of vertex~$j$ on vertex~$i$, and \adbb{is} assumed to take the value zero if $i$ and~$j$ are
not neighbours, i.e.\ if $A_{ij} = 0$; \adbb{we also assume that $W_{ii} = 0$ for each~$i$.}  
The parameter~$\rho$ can be varied, and $|\rho|$ mirrors the overall strength of the 
\adbb{dependence  between the observations at neighbouring vertices.} 

The network disturbance model can be estimated by using the method of maximum likelihood (MLE);  details of the
procedure are given in Section~\ref{Details}.
\adb{In classical settings, including repeated independent observations from a common underlying distribution,
the maximum likelihood estimator is asymptotically unbiased, as the number of observations tends to infinity. }
Here, the observations are dependent, and the standard asymptotic conclusions need no longer be valid.
Indeed, it has been widely observed in simulation studies that, as the underlying network becomes denser, 
and the corresponding weight matrix~$W$ becomes less sparse,
the maximum likelihood estimator~$\rhat$ of~$\r$ becomes more variable and more negatively biased. This has
been extensively documented, \adbb{largely in the context of the related network effects models
of \cite{doreian1989network},} \adbb{defined in~\Ref{NEffects-model} below,}
 for instance in \cite{mizruchi2008effect},  \cite{smith2009estimation}, \cite{farber2010topology},  
\cite{neuman2010structure} and \cite{la2018finite}.  \adbb{Methods for correcting the bias of the MLE, using asymptotics
appropriate in settings where consistent estimation is possible, have since been proposed in} 
\adbb{\cite{bao}, \cite{yubaiding} and \cite{yang}.  Other methods of estimation have been developed by 
\adb{\cite{KP1999}, who introduce a generalized method of moments estimator,} by \cite{dittrich}, with
a Bayesian approach, and by \cite{strumann}, whose techniques are based on Hodges--Lehmann estimators.}
In \cite{smith2009estimation}, theoretical investigation 
shows that, at least for networks very close to the complete graph, the MLE may \adbb{lead to misleading results:
see the final example in Section~\ref{Ill-conditioned} for a fuller explanation,} \adbb{that makes clear why
no estimator of~$\r$ can work well in such circumstances.}

In this paper, we examine the problem in greater depth, showing why the MLE can be expected to be biased and 
lacking in precision, whenever the underlying network is dense.
In Section~\ref{Cramer-Rao}, we use Cram\'er--Rao theory to show that, if the underlying network is 
sufficiently dense, then {\it any\/} (more or less unbiased)
estimator of~$\r$ must necessarily be unreliable. 
The bound is expressed as a lower bound~\Ref{CR} on the variance of any (unbiased) estimator, 
and is given explicitly as a function
of the matrix~$W$ of weights and of the true value of~$\rho$.  The structure of the network enters by way of~$W$.
For instance, the choice $W_{ij} := A_{ij}/\deg(i)$ makes explicit the dependence of the weights on the adjacency 
matrix~$A$ of the network.
We then, in Section~\ref{MLE-section}, examine the behaviour of the maximum likelihood estimator in enough
detail to explain the empirical observations.  In particular, we are able to show that the MLE of~$\rho$ is biased
when the network is sufficiently dense.  \red{In ordinary least squares estimation, the estimator of~$\b$ has
small variance if the smallest eigenvalue of the matrix~$X\TT X$ is large.  Here, it is possible, for rather special 
graphs, that the MLE of~$\b$ may still have substantial variance, even when this condition is satisfied;  see Examples
1 and~2 in Section~\ref{Ill-conditioned}.}

The exact formulae for the bias and variance of the MLE are
in general indigestible, but we discuss some special cases in which it is possible to estimate their
magnitudes. 
In Section~\ref{Quadratic-form-estimation}, we consider a
simpler and more intuitive estimator of~$\r$, that has better properties as regards bias, but is still
lacking in precision in dense graphs, as it must be, in view of the lower bound
given in Section~\ref{Cramer-Rao}.  Its properties are illustrated by some simulations.
\red{An R~function for computing the estimator, as well an empirical measure of its precision,
is available on request from the authors.}

There is a related model, the {\it Network Effects Model\/} of \cite{doreian1989network}
\adbb{(see also \cite[(4.1)]{ord1975}}, which takes the form
\eq\label{NEffects-model}
      K_\r y \Eq X\b + \nu,
\en
where the components of~$\nu$ are realizations of independent normal random variables with mean zero and
variance~$\s^2$, \adbb{and~$K_\r$ is as defined in~\Ref{NDM-errors-2}}.  This differs from the network 
disturbance model because, rewriting it as
\eq\label{NEffects-model-1}
    y \Eq K_\r^{-1}X\b + \e,
\en
the errors in~$y$ have the same dependence as before, but the {\it model structure\/} explaining~$y$ also
changes with~$\r$.  The limitations as regards estimating~$\rho$ in dense networks are analogous, and for
\adbb{similar reasons. We discuss the model briefly in \adbb{Section~\ref{NE-models-sect}.}}

\red{The main conclusions of the paper are sketched in Section~\ref{conclusion}, and additional calculations
can be found in the Appendix. }

\section{Details of the model}\label{Details}
\setcounter{equation}{0}

The maximum likelihood estimates of the parameters \adbb{in the simple linear regression model
$\by \Eq  X \b + \beps$ of~\Ref{linear-regression}, with independent and identically distributed errors,}  
are given by
\eq\label{OLS-MLE}
    \hat\b \Eq \hat\b(y) \Eq (X\TT X)^{-1}X\TT y;\qquad \hat\s^2 \Eq \hat\s^2(y) \Eq n^{-1} y\TT (I - H(X)) y,
\en
where 
\eq\label{H-proj-def}
    H(X) \Def X(X\TT X)^{-1}X\TT
\en 
is symmetric and idempotent, and represents the orthogonal projection
from~$\re^n$ onto the linear span of the columns of~$X$. 
 Writing 
\eq\label{OLS-model}
    M_0\colon\ \by \Eq X\b + \beps
\en
to denote the probabilistic model from which
the data are sampled, the random quantities $\hat\b(\by)$ and~$\hat\s^2(\by)$ have distributions given by
\eq\label{OLS-MLE-2}
       \hat\b(\by) \ \sim\ \MVN_m(\b,\s^2(X\TT X)^{-1});\qquad  \hat\s^2(\by) \ \sim\ n^{-1}\s^2 \chi^2_{n-m}, 
\en
\adbb{where $\MVN_m(\mu,\Sigma)$ denotes the $m$-dimensional multivariate normal distribution with mean~$\mu$
and covariance marix~$\Sigma$, and~$\chi^2_l$ denotes the chi-squared distribution with~$l$ degrees of
freedom. The two estimators are also} independent of one another \adb{(see \cite[Section~4.b]{rao}).}

\adbb{In the network disturbance model $\by \Eq  X \b + \beps$, the errors~$\beps$ are given by~\Ref{NDM-errors-2},
with dependence deriving from the weight matrix~$W$.}
 It is often convenient to assume that the row sums of~$W$ are all equal to~$1$, in which case the quantity
$\sum_{j=1}^n W_{ij}\beps_j$ represents an average of the errors at the neighbours of vertex~$i$;  then~$\rho$
would typically be smaller than~$1$ in modulus, so as to ensure that the average influence of the neighbours is not
larger than the size of a typical error.  A simple choice for~$W_{ij}$ is $A_{ij}/\deg(i)$, where~$\deg(i)$ denotes
the degree of the vertex~$i$.  \adb{If there are isolated vertices~$i$, the corresponding row sums cannot be made
equal to~$1$, since then $W_{ij} = 0$ for all $j\ne i$.}

The implicit representation~\Ref{NDM-errors-2} of the errors can be written as a direct formula for~$\beps$
in terms of~$\bnu$, as
\eq\label{eps-from-nu}
       \beps \Eq K_\r^{-1}\bnu,
\en
provided that the matrix~$K_\r$ can be inverted.  The expansion
\eq\label{eps-from-nu-expansion}
    K_\r^{-1} \Eq (I - \rho W)^{-1} \Eq I + \sum_{l\ge1}\r^l W^l
\en
makes sense provided that the sum $\sum_{l\ge1}\r^l W^l x$ converges absolutely for each $x \in \re^n$;
\adbb{that is, if~$\adb{|\r|}$ is less than~$1/r(W)$, where~$r(W)$ is the spectral radius of~$W$.} 
In particular, if~$W$ has all of its row sums equal to~$1$, the vector~$\bone$ with all components equal to~$1$
is a right eigenvector with eigenvalue~$1$, implying that $r(W)=1$ 
\adbb{(see \cite[Corollary 8.1.30, p.493]{horn2012matrix})}.  \adbb{Then}
\[
     \beps \Eq K_\r^{-1}\bnu \Eq \bnu + \sum_{l\ge1}\r^l W^l\bnu,\qquad |\r| < 1,
\]
represents~$\beps_i$ as a weighted sum of the contributions of the different elements~$\bnu_j$,
averaged along paths in the network from $j$ to~$i$, and with overall weight~$\r^l$ for paths of length~$l$
(possibly negative, if $\r < 0$).  If $|\r| > 1/r(W)$, this interpretation no longer holds, since the
expansion~\Ref{eps-from-nu-expansion} is no longer well defined. The resulting model is then not
so easy to understand, since, at least locally, the cumulative effect on the observation at~$i$ from those at
distance~$d$ from~$i$ increases geometrically with~$d$.  In time series analysis, the analogue is that of processes 
that are `explosive', and the analogue of the condition $|\r| < 1/r(W)$ is \adb{used to ensure that} a time series is
stationary.  For this reason, 
we shall assume in what follows that $|\r| < 1/r(W)$, so
that the inverse of~$K_\r$ exists, and is explicitly given by~\Ref{eps-from-nu-expansion}.

The model $\by \Eq  X \b + \beps$ with errors given by~\Ref{NDM-errors} is simple to analyze if the value of~$\rho$
is known, since it can be written in the form
\eq\label{M-rho-def}
     M_\r\colon \quad K_\r \by \Eq K_\r X \b + \bnu,
\en
with the components of~$\bnu$ being independent normal~$\cN(0,\s^2)$ random variables, and can thus be fitted as
a standard linear regression based on data~$y' := K_\r y$.  Equations~\Ref{OLS-MLE} then give the estimates
\eqa
    \hat\b_\r &=& \hat\b_\r(y) \Eq ((K_\r X)\TT K_\r X)^{-1}(K_\r X)\TT K_\r y; \label{NDM-MLE}\\
          \shat^2_\r &=& \shat^2_\r(y) \Eq \adbb{n^{-1}\{K_\r(y - X\hat\b_\r)\}\TT \{K_\r(y - X\hat\b_\r)\}}
       \Eq n^{-1} (K_\r y)\TT (I - H_\r) K_\r y, \non
\ena
where, \adb{with~$H$ as in~\Ref{H-proj-def},} \adbb{
\eq\label{H-rho-def}
    H_\r \Def H(K_\r X) 
\en 
}
\adbb{and so $K_\r X \hat\b_\r = H_\r K_\r y$. The residuals after fitting~$M_\r$ are thus
\eq\label{residuals}
     \hbe_\r(y) \Def y - X\bhat_\r(y) \Eq K_\r^{-1}(I - H_\r)K_\r y.
\en
}
The distributions of $\hat\b_\r(\by)$ and~$\hat\s^2_\r(\by)$ under the model~$M_\r$ are given by
\eq\label{NDM-MLE-2}
       \bhat_\r(\by) \ \sim\ \MVN_m(\b,\s^2((K_\r X)\TT K_\r X)^{-1});
                       \quad  \shat^2_\r(\by) \ \sim\ n^{-1}\s^2 \chi^2_{n-m}, 
\en
and the two estimators are independent of one another. \adbb{In particular, both $\bhat$ and~$\shat^2$
have precision much as in simple linear regression, and (under
reasonable assumptions about~$X$) can be expected to become more accurate as~$n$ increases.
If the true model is~$M_{\r'}$, for $\r' \neq \r$, the expectation of~$\bhat_\r(\by)$ is still equal to~$\b$,
and, although the variances of its components may be larger, they are typically comparable to those
on the true model,} \adb{again under reasonable conditions;  see Appendix~\ref{off-the-model} for more detail.}

If~$\r$ is not known, the estimates~\Ref{NDM-MLE} still maximize
\adb{the log-likelihood 
\eq\label{LL}
  \ell(\th;y) \Def -\tfrac n2 \{\log(\s^2) + \log (2\pi)\} - \frac1{2\s^2} \{K_\r(y-X\b)\}\TT  \{K_\r (y-X\b)\} 
               + \log\det K_\r 
\en
of the model~$M_\r$ of~\Ref{M-rho-def}, for each given value of~$\r$,
where $\th := (\b,\s^2,\r)$ denotes the full set of parameters; see \cite[Equation~(1.3)]{doreian1989network}.
The resulting maximum of the log-likelihood, for fixed~$\r$,} is given by
\eq\label{NDM-ML-max}
     -\tfrac n2 \{\log(\shat^2_\r(y)) + \log (2\pi) + 1\} + \log\det K_\r;
\en
\adbb{this can be deduced by observing that, with~$\b$ replaced by~$\hat\b_\r$, the middle term of~\Ref{LL} 
reduces to~$-n/2$, in view of~\Ref{NDM-MLE}.}
Indeed, the first term in~\Ref{NDM-ML-max} is the maximum of the likelihood for ordinary least squares applied 
to the data~$y' := K_\r y$,
and the determinant enters because the probability densities $f_{\by'}$ of~$\by'$ and~$f_{\by}$ of~$\by$ are related by 
$f_{\by}(v) = \det K_\r\, f_{\by'}(K_\r v)$.  Then~$\r$ can be chosen to (numerically) maximize the expression 
in~\Ref{NDM-ML-max}, or, equivalently, to minimize
\eq\label{infer-disturbance-rho-hat}
   F(\r;y) \Def \log(\hat\s^2_\r(y)) - \tfrac2n \log\det K_\r, \qquad |\r| < 1/r(W).
\en

\section{A general limit to precision}\label{Cramer-Rao}
\setcounter{equation}{0}

The $\rho$-component of the score function is obtained by differentiating the log-likelihood~\Ref{LL} with respect to~$\r$.
Noting that $\frac{d}{d\r}K_\r = - W$, and that,
if the eigenvalues of~$W$ are denoted by $\l_1,\ldots,\l_n$, we have $\det K_\r = \prod_{i=1}^n (1 - \r\l_i)$.
This gives the $\rho$-component of the score function as
\eqa
   s\urho(\th;y) &:=& \frac{\partial}{\partial\r}\ell(\th;y)
           \Eq \frac1{\s^2} (y-X\b)\TT K_\r\TT W (y-X\b) - \sum_{l=1}^n \frac{\l_l}{1-\r\l_l}  \non\\
      &=& \frac1{\s^2} \bigl(\{K_\r(y-X\b)\}\TT  Z_\r \{K_\r (y-X\b)\}\bigr) - \Tr\{Z_\r\}, \label{urho-def}
\ena
where we define the matrix \adbb{ 
\eq\label{Z-rho-def}
   Z_\r \Def WK_\r^{-1}. 
\en
} 
Now, if~$C$ is any $n\times n$ matrix, and~$\tbnu$ denotes an $n$-vector of independent
{\it standard\/} normal random variables, then the quadratic form $\tbnu\TT C \tbnu$ has mean 
\eq\label{quad-form-mean}
   \ex\{\tbnu\TT C \tbnu\} \Eq \Tr\{C\}.
\en
Furthermore, for {\it symmetric\/}~$C$,
since $\cov\{\tbnu_i\tbnu_j,\tbnu_k\tbnu_l\} = 0$ unless either $i \ne j$ and $\{i,j\} = \{k,l\}$ or $i=j=k=l$,
the variance of $\tbnu\TT C \tbnu$ is given by
\eqa
    \var\{\tbnu\TT C \tbnu\} &=& \sum_{i=1}^n \sum_{j\ne i} C_{ij}(C_{ij} + C_{ji}) + \sn C_{ii}^2(\ex\{\tbnu_1^4\} - 1)
                                          \label{quad-form-1} \\
                           &=& 2 \Tr(C^2),   \label{quad-form-var}
\ena
because $\ex\{\tbnu_1^4\} = 3$ when $\tbnu_1 \sim \cN(0,1)$; see \cite[Chapter~2, Theorem 4 (i) \& (ii)]{searle1971linear}.
If~$C$ is not symmetric, the 
expression~\Ref{quad-form-var} can be used with~$C$ replaced by the symmetric matrix $C' := \half(C + C\TT)$, 
since $\tbnu\TT C' \tbnu = \tbnu\TT C \tbnu$.
Hence, taking expectations on the model with parameters $\th_0 := (\b_0,\s^2_0,\r_0)$, we have 
\[
    \ex_{\th_0}\, s\urho(\th_0;\by) \Eq \Tr\{Z_{\r_0}\} - \Tr\{Z_{\r_0}\}\ \Eq 0,
\]
as it must, and
\eq\label{var-score}
   \var_{\th_0}\, s\urho(\th_0;\by) \Eq 2 \Tr\bigl\{\bigl( \half(Z_{\r_0} + Z_{\r_0}\TT) \bigr)^2\bigr\} 
                \Eq \Tr\{Z_{\r_0}^2 + Z_{\r_0} Z_{\r_0}\TT\}.
\en
\adbb{This, together with \cite[Theorem~6.2]{hardle2015applied} and \cite[Theorem~6.6 and (6.25)]{LehmannCasella},
implies the following limit on the precision of estimators of~$\r$.} 

\bigskip
\nin{\bf Cram\'er--Rao lower bound.} \adbb{The variance of any unbiased estimator of~$\rho$ is no smaller than
\eq\label{CR}
     \g_{\r_0}^2 \Def \frac1{\Tr\{Z_{\r_0}^2 + Z_{\r_0} Z_{\r_0}\TT\}}\,,
\en
where $Z_\r := W K_\r^{-1}$ and~$K_\r$ is as in~\Ref{eps-from-nu}.
}

\bigskip
\nin In practice, the bound \adbb{implies a precision that cannot be significantly improved upon for any estimator of~$\rho$,
without introducing substantial bias;}
even the best `sensible' estimator of~$\rho$ has to have random variability on the scale of~$\g_{\r_0}$.

Note that, for any $n\times n$ matrix~$M$, 
\eq\label{Cauchy-Schwarz}
     |\Tr\{M^2\}| \Eq \Bigl|\sn\sjn M_{ij}M_{ji}\Bigr| \Le \sn\sjn M_{ij}^2 \Eq \Tr\{MM\TT\},
\en
by Cauchy--Schwarz, so that the two traces appearing in~$\g_{\r_0}^2$ are not in general the same,
and that
\eq\label{CR-element}
   \Tr\{M^2 + MM\TT\} \Eq \sn\sjn \{M_{ij}M_{ji} + M_{ij}^2\} \Eq \half\sn\sjn (M_{ij} + M_{ji})^2 \ \ge\ 0.
\en

\paragraph{Consequences of~\Ref{CR}: the case $\r_0=0$.}
As a simple example of what is to be expected, suppose that $\r_0 = 0$, and that the underlying graph
has no isolated vertices, so that $\deg(i) \ge 1$ for all~$i$.  Let $W := D^{-1}A$,
\adbb{where~$D$ is the diagonal matrix with elements $D_{ii} := \deg(i)$, $1\le i\le n$.} Then $Z_{\r_0} = Z_0 = W$, 
and we can be specific about the traces in~\Ref{CR}:
\eqa
    \Tr\{W^2\} &=& \sn\sjn \frac{A_{ij}}{\deg(i)\deg(j)} \non\\
    &\le& 
           \Tr\{WW\TT\} \Eq \sn\sjn \frac{A_{ij}}{\deg(i)^2} \Eq \sn\frac1{\deg(i)} \Le n.  \label{W-traces}
\ena
It follows from the equality in~\Ref{W-traces}
that $\Tr\{WW\TT\} \ge n/\bard$, where~$\bard  := n^{-1}\sn \deg(i)$ denotes the average degree of 
a vertex, and hence
that $\g_{0}^2 \le \bard/n$.  However, if the vertex degrees are reasonably homogeneous, in the sense all the vertices have 
degrees at least~$c\bard$, for some $0 < c < 1$, then it follows from~\Ref{W-traces} that $\Tr\{WW\TT\} \le n/c\bard$,
and thus, from~\Ref{CR}, that $\g_{0}^2 \ge c\bard/n$. Hence, in these simple circumstances, the variability of any
sensible estimator of~$\rho$ is on a scale of at least~$(\bard/n)^{1/2}$.  

\adb{
These calculations have particular implications for dense networks.  If $W = D^{-1}A$, and if
all vertices have degrees exceeding~$cn$, for some fixed $c > 0$, it follows that $\g_{0} \ge \sqrt c$, and it 
is not possible to obtain an accurate estimate of~$\r$ by any method, no matter how large
the size~$n$ of the network, if in fact $\r_0=0$.  For more general choices of~$W$, 
$\Tr\{WW\TT\} = \sum_{i,j=1}^n W_{ij}^2$, and
if all the elements of~$W$ are smaller than~$C/n$, for some~$C < \infty$, it follows that $\Tr\{WW\TT\} \le C^2$,
and hence, from~\Ref{CR}, that \red{$\g_{0}^2 > 1/C^2$.}  Hence, here too, any estimator of~$\rho$ has 
significant variability on the scale of~$1/C$, whatever the size~$n$ of
the network, if in fact $\r_0 = 0$.
}

\paragraph{Consequences of~\Ref{CR}: general $\r_0$.}
If $\r_0 \ne 0$, analogous results are true.
Supposing that  $\sjn W_{ij} = 1$ for each~$i$, it follows that, for all $l \ge 1$, $\sjn (W^l)_{ij} = 1$ for all~$i$ also.
Hence $(W^l)_{ij} = \sum_{k=1}^n (W^{l-1})_{ik}W_{kj}$ is an average of the values $W_{1j},\ldots,W_{nj}$,
Thus, if all the elements of~$W$ are smaller than~$C/n$, it follows that 
$(W^l)_{ij} \le C/n$ also.  Hence it is immediate that, for $l,l' \ge 1$, 
\eq\label{W-trace-powers-1}
     \Tr\{W^l\} \Le C \quad\mbox{and}\quad \Tr\{W^l(W\TT)^{l'}\} \Le C.
\en
Thus, for $\r \ge 0$, we have
\eqa  
      \Tr\{Z_\r^2\} &\le& \Tr\{Z_\r Z_\r\TT\} \Eq \Tr\{W(I-\r W)^{-1}(I-\r W\TT)^{-1}W\TT\} \non\\
              &\le& C\sum_{l\ge0}(l+1)\r^l \Eq \frac C{(1-\r)^2}\,, \label{W-trace-powers-1.5}
\ena
and so
\[
     \g_{\r_0} \ \ge\ C^{-1/2}(1-\r) .
\]
Once again, there is substantial variability in the estimate of~$\r$ unless $C := n\max_{i,j}W_{ij}$ is large. 
If~$\r < 0$, the bounds in~\Ref{W-trace-powers-1} can be used to justify similar conclusions, by bounding
the difference $\Tr\{Z_\r Z_\r\TT\} - \Tr\{WW\TT\}$, provided that~$|\r|$ is not too large.

\adb{
In particular, in uniformly dense networks, if $W = D^{-1}A$, there is substantial variability in any estimator
of~$\r$ if all vertex degrees exceed~$n/C$, and~$C$ is not large.}

\section{Maximum likelihood estimation}\label{MLE-section}
\setcounter{equation}{0}

\adb{
The discussion in the previous section shows that there is an upper limit on the precision of any
(unbiased) estimator of~$\rho$, and that, in uniformly dense networks, any such estimator must have
substantial variability. In this section, we go into more detail about the properties of the MLE.
}

\subsection{General theory}

\nin{\bf An equation for the MLE} 

\msk
To minimize the quantity~$F(\r;y)$ of~\Ref{infer-disturbance-rho-hat}, it is usual to look for values of~$\r$ 
such that 
\eq\label{U-rho-eqn}
  U(\r;y) \Def \frac{d}{d\r} F(\r;y) \Eq 0.  
\en
Since~$F(\r;y)$ is explicitly given as a function of~$\r$, its derivative can be calculated, giving
\eqa
 U(\r;y) &=& \frac1{\shat^2_\r(y)}\,\frac{d\shat^2_\r(y)}{d\r}  - \frac2n\, \frac d{d\r}\,\log\det K_\r 
               \non\\
    &=& -\frac1{n\shat^2_\r(y)}\, (K_\r y)\TT (I - H_\r)(Z_\r\TT + Z_\r)(I-H_\r) K_\r y  
           + \frac2n \sum_{l \ge 1}\r^l\Tr\{W^{l+1}\}; \phantom{XX} \label{D-F}
\ena
the details are given in~Appendix~\ref{U-calc}.  
\adb{The maximum likelihood estimate~$\hrho\ML$ of~$\r$ is now typically obtained by solving $U(\r;y) = 0$;}
\red{Example~3 in Section~\ref{Ill-conditioned} shows that this may not always work.}

\adbb{Note that the data enter through a ratio of quadratic forms
in~$y$.  Under~$M_{\r_0}$, as observed in~\Ref{NDM-MLE-2}, $\shat^2_{\r_0}(\by) \ \sim\ n^{-1}\s^2 \chi^2_{n-m}$, 
having mean $\s^2(1 - m/n)$ and variance less than~$2n^{-1}\s^4$, so that, to a reasonable approximation,
$\shat^2_{\r_0}(y)$ can be replaced by~$\s^2$ in~\Ref{D-F} if~$M_{\r_0}$ is true, provided that~$n$ is large.
As a result, $U(\r;y)$ depends on the data~$y$ largely through a single quadratic form in~$y$.}

\bsk
\nin{\bf \adbb{Approximating the} distribution of $U(\r;\by)$ \adb{under $M_{\r_0}$}}

\msk
In order to gauge the typical values of $U(\r;y)$ to be obtained under the true model~$M_{\r_0}$,
we express~$U(\r;\by)$ in terms of the (not directly observable) errors~$\bnu$, rather than in terms of~$\by$.
Under~$M_{\r_0}$, $\by = X \b + K_{\r_0}^{-1}\bnu$, and it thus follows from \Ref{infer-disturbance-rho-hat} 
and~\Ref{D-F} that
\eqa
    U(\r;\by) &=& \frac1{\shat^2_\r(\by)}\,\frac{d\shat^2_\r(\by)}{d\r}  - \frac2n\, \frac d{d\r}\,\log\det K_\r \non\\
    &=& -\frac1{n\shat^2_\r(K_{\r_0}^{-1}\bnu)}\, \bnu\TT Q(\r;\r_0)\bnu +  \frac2n \sum_{l \ge 1}\r^l\Tr\{W^{l+1}\}, 
            \label{D-F-nu}
\ena
where
\eq\label{Q-def}
   Q(\r;\r_0) \Def (K_\r K_{\r_0}^{-1})\TT (I - H_\r)(Z_\r\TT + Z_\r)(I-H_\r) K_\r K_{\r_0}^{-1}.
\en
Note that~$\b$ does not appear in~\Ref{D-F-nu} because $(I-H_\r)^T = I - H_\r$ and $(I-H_\r)K_\r X = 0$,
\adbb{implying in particular, using~\Ref{NDM-MLE}, that $\shat^2_\r(\by - X\b) = \shat^2_\r(\by)$.}  
Hence, \adb{and because 
$\shat^2_\r(K_{\r_0}^{-1}\bnu)$ is close to~$\s^2$ for large~$n$,} 
 the random part of $U(\r;\by)$ is given, to a close approximation, by $ -\bnu\TT Q(\r;\r_0)\bnu/(n\s^2)$.

\adbb{In the remaining argument, we make frequent use of elementary properties of the trace, that $\Tr(AB) = \Tr(BA)$
and that $\Tr(A\TT) = \Tr(A)$, as well as noting that \adbb{$z\TT M z = z\TT M\TT z$ 
for any compatible $z$ and~$M$,} and that $H_\r$ and~$I - H_\r$
are idempotent.}
 The expectation of $\bnu\TT Q(\r_0,\r_0) \bnu$, using~\Ref{quad-form-mean}, is given by
\eqa
    \ex\{\bnu\TT Q(\r_0,\r_0) \bnu\} &=& \s^2 \Tr\{Q(\r_0,\r_0)\} 
                   \Eq \adbb{2\s^2 \Tr\{(I- H_{\r_0})Z_{\r_0}(I- H_{\r_0})\}} \non\\
      &=& 2\s^2 \Tr\{(I- H_{\r_0})Z_{\r_0}\} \Eq 2\s^2\bigl(\Tr\{Z_{\r_0}\} - \Tr\{H_{\r_0}Z_{\r_0}\}\bigr), \label{Q-nu-mean}
\ena
\adbb{because $K_{\r_0} K_{\r_0}^{-1} = I$, and using $\Tr(AB) = \Tr(BA)$ and $(I-H_{\r_0})^2 = I-H_{\r_0}$.}
Then, using $\Tr \{W\} = 0$,
\eq\label{trace-Z_rho}
    \Tr\{Z_{\r_0}\} \Eq \Tr\{W(I - \r_0 W)^{-1}\} \Eq \sum_{l \ge 1}\r_0^l \Tr\{W^{l+1}\}. 
\en
Thus we have
\[
   -\ex\bigl\{\s^{-2}\bnu\TT Q(\r_0;\r_0)\bnu\bigr\} +  2\sum_{l \ge 1}\r^l\Tr\{W^{l+1}\}
             \Eq 2\Tr\{H_{\r_0} Z_{\r_0}\}.
\]
Hence, from \adbb{\Ref{D-F-nu} and~\Ref{Q-nu-mean}}, neglecting the small variability in $\shat^2_\r(\by)$, 
 the quantity $n U(\r_0;\by)$ under the model~$M_{\r_0}$ satisfies
\eq\label{mu0-def}
    \ex_{\r_0}\{n U(\r_0;\by)\} \ \approx\ \m_0 \Def \m(\r_0) \Def 2\Tr\{H_{\r_0} Z_{\r_0}\}.
\en

The principal source of variability in $n U(\r_0;\by)$ arises from 
$\s^{-2}\bnu\TT Q(\r_0;\r_0)\bnu$, whose variance under~$M_{\r_0}$, from~\Ref{quad-form-var}, is given 
by $2\Tr\{Q(\r_0;\r_0)^2\}$.   This yields
\eqa
     \var_{\r_0}(n U(\r_0;\by)) &\approx& 2\Tr\{Q(\r_0;\r_0)^2\} \non\\
      &=& \adbb{2 \Tr\{(Z_{\r_0} + Z_{\r_0}\TT) (I - H_{\r_0}) (Z_{\r_0} + Z_{\r_0}\TT) (I - H_{\r_0})\}} \non\\
      &=:& \t_0^2 \Def \t_0^2(\r_0). 
                      \label{tau0-def}
\ena
Combining this with~\Ref{mu0-def}, \adbb{we reach the following}

\msk\nin\adbb{{\bf  Conclusion~1:} the  distribution of~$nU(\r_0;\by)$ can be written in the form}
\eq\label{nU-approx-values}
        nU(\r_0;\by) \Eq \m_0 + \t_0 \EEE(\by),
\en
where $\EEE(\by)$ has mean approximately~$0$ and variance approximately~$1$.

\bsk
\nin{\bf Consequences for the distribution of the MLE}

\msk
\adb{
Using the information above, we can now make approximate deductions about the MLE $\hrho := \hrho\ML$, obtained by setting
$U(\r;y)$ equal to zero.  We begin by noting} that $U(\hrho,y) = 0$, whereas, as in~\Ref{nU-approx-values},
$nU(\r_0;y) = \m_0 + \t_0 \EEE(y)$, where~$\EEE(\by)$ approximately has mean zero and variance~$1$.  
Writing~$U'(\r;y) := \frac{d}{d\r} U(\r;y)$, this implies that
\[
     -n^{-1}(\m_0 + \t_0 \EEE(y)) \Eq \{U(\hrho;y) - U(\r_0;y)\} \Eq (\hrho-\r_0)U'(\r_0;y) + o(|\hrho - \r_0|),
\]
\adb{where the final error term is relatively small}
if $(|\m_0| + \t_0)/n|U'(\r_0;y)|$ is small enough, leading to the approximate expression
for the error $\hrho - \r_0$ in estimating~$\r_0$:
\eq\label{rhohat-error}
     \hrho\ML -\r_0 \ \approx\ -\frac{\m_0 + \t_0 \EEE(y)}{nU'(\r_0;y)}\,.
\en
\adbb{
This is essentially a one--step \adbb{Newton--Raphson} approximation to the value of~$\r$ that 
solves $U(\r;y) = 0$, starting from the guess~$\r_0$.}
Thus it is necessary also to examine the magnitude of~$nU'(\r_0;y)$.  Since~$U(\r;y)$ is explicitly given,
there are exact formulae for its derivative, given in Appendix~\ref{U-deriv};  however, they are not easy to understand.
Instead, it is possible to approximate the {\it expectation\/} of~$nU'(\r_0;\by)$ on the model~$M_{\r_0}$, giving
\eq\label{nU-dash-approx}
     n\ex_{\r_0} U'(\r_0;\by) \ \approx\ \D_0,
\en
where
\eqa
  \D_0  &=& \half\t_0^2 + 2\Tr\{ - Z_{\r_0} H_{\r_0}(Z_{\r_0} + Z_{\r_0}\TT)(I - H_{\r_0})  +  H_{\r_0} Z_{\r_0}^2\}; 
                    \label{U-dash-expec}                                             
\ena
the details are again to be found in Appendix~\ref{U-deriv}.
Thus, from~\Ref{rhohat-error}, we have the following

\msk\nin{\bf Conclusion~2:} the difference between $\hrho\ML$ and~$\r_0$ can be expected to be of
magnitude roughly
\eq\label{rhohat-error-2}
             \frac{|\m_0| \adb{\pm} \t_0}{\D_0}\,,
\en
with the term $|\m_0|/\D_0$ corresponding to systematic bias, and $\ssc_0 := \t_0/\D_0$ reflecting the scale
of the randomness in the estimator.  

\msk
\adbb{Note that the measure $\ssc_0$ should be treated only as a rough
guide to the actual variability of~$\hrho\ML$ under~$M_{\r_0}$, since~\Ref{rhohat-error} is only a linear
approximation to the solution of a non-linear equation, and because there may be added variability
arising from the variation of~$nU'(\r_0;\by)$ about its mean.}

\subsection{Estimating~$\r$ when there is no \adb{structural element}}\label{no-structure}
To get an idea of the relative magnitudes of the terms $\m_0$, $\tau_0$ and~$\D_0$
in~\Ref{rhohat-error-2}, it is helpful to look
at some special cases.  For the first, suppose that there is no structural element at all --- not even a mean to
be fitted --- \adb{which corresponds to taking $X=0$; in this case, the formulae \Ref{mu0-def}, \Ref{tau0-def} 
and~\Ref{U-dash-expec} should be interpreted as having $H_\r=0$, the projection onto the space~$\{0\}$.} 
Then $\m_0 = 0$,  so that any bias is small,
\[
     \t_0^2 \Eq 4 \Tr\bigl\{ Z_{\r_0}^2 + Z_{\r_0}Z_{\r_0}\TT \bigr\},\quad\mbox{where}\quad
                Z_{\r} \Eq W(I-\r W)^{-1},
\]
and
\[
    \D_0 \Eq \t_0^2/2,
\]
giving
\eq\label{t0-D0}
   \ssc_0 \Eq \frac{\t_0}{\D_0} \Eq \frac1{\sqrt{\Tr\{Z_{\r_0}^2 + Z_{\r_0}Z_{\r_0}\TT\}}}.
\en
This is the same as~$\g_{\r_0}$ in~\Ref{CR}.
Thus, from~\Ref{rhohat-error-2}, the typical magnitude of $|\hrho\ML - \r_0|$ is of order~$\g_{\r_0}$.
This is the best order of error that could theoretically have been hoped for.  Nonetheless, as observed above, it
may not lead to accurate estimation of~$\rho$, \adb{for instance if $W = D^{-1}A$ and the underlying network 
is sufficiently dense, even though the {\it bias\/} is small.}

\subsection{Estimating~$\r$ when structure is present}\label{Estn-with-structure}

To accommodate the effect of having structure in the model, suppose first that only an overall mean is to be fitted,
so that~$X$ is the $n$-vector $\bone$, having all elements equal to~$1$.  Suppose also that there are no isolated vertices, 
and that~$W$ is chosen to have $\sjn W_{ij} = 1$ for all~$j$, so that $W\bone = \bone$ \adb{and the spectral radius $r(W)=1$}.
\red{Then it is immediately verified that $H_{\r} = H_* := n^{-1}\bone\bone\TT$ for all~$\r$,
because, for instance, $K_\r X = (I-\r W)\bone = (1-\r)\bone$.  Hence $WH_\r = H_*$ and $Z_\r H_{\r} =  (1-\r)^{-1}H_*$. 
As a result, by arguments that are made precise in Appendix~\ref{4.3-extras},} \Ref{tau0-def} and~\Ref{U-dash-expec} yield
\eqa
          \t_0^2  &=&  4\Tr\{(Z_{\r_0}Z_{\r_0}\TT + Z_{\r_0}^2)(I - H_*) \} 
                                \adb{ + 2\{(1-\r_0)^{-2} - n^{-1}|Z_{\r_0}\TT \bone|^2\} } \non\\
                  &\le&  \adb{ 4\Tr\{(Z_{\r_0}Z_{\r_0}\TT + Z_{\r_0}^2)(I - H_*) \} }; \label{with-structure} \\
          \D_0  &=&  \half\t_0^2  +  2(1-\r_0)^{-2} ; \non
\ena
\red{here, and in what follows,  $|x|$ for a vector~$x$ denotes the Euclidean norm.}
Hence, \adbb{if $\r_0$ is fixed and less than~$1$, $\ssc_0$ is of order $1/(\t_0 + 1)$,
and this, from~\Ref{rhohat-error-2}, represents the scale of the {\it random\/} variation in the estimate of~$\r_0$. 
Now $1/(\t_0 + 1)$} is small if~$\t_0^2$ is large compared to~$1$,
but can be significantly larger than the lower bound~$\g_{\r_0}$ in~\Ref{CR}, if $\t_0^2$ is much smaller than 
$1/\g_{\r_0}^2 = \Tr\{Z_{\r_0}Z_{\r_0}\TT + Z_{\r_0}^2\}$.
As in the case where there is no structure, \adb{if $W = D^{-1}A$ and the graph is sufficiently dense,} then
$\g_{\r_0}$ need not be small, in which case the MLE has substantial variability.

In contrast to the situation without structure, the principal element~$\m_0$ of the bias in the estimation of~$\r$
need not be zero.  If $X = \bone$,
it follows \adbb{from~\Ref{mu0-def}} that
\[
    \m_0 \Eq 2\Tr\{Z_{\r_0}H_{0}\} \Eq 2(1-\r_0)^{-1} \ >\ 0.
\]
\adb{In view of Conclusion~2,
the bias becomes important when~$\D_0$ is not large, which, from~\Ref{with-structure}, is when~$\t_0^2$ is not large;
as discussed above, this occurs \red{for example} if $W = D^{-1}A$ and the graph is uniformly dense.}
Note also that, from \Ref{rhohat-error} and~\Ref{rhohat-error-2}, the value of the difference~$\hrho\ML -\r_0$ resulting 
from~$\m_0$ is roughly $- 2/\{\D_0(1-\r_0)\}$, which typically implies a negative bias, as demonstrated 
in the results of the simulations referred to in the introduction.

If there is more structure in the regression, but still assuming that \adb{the first column of~$X$ is~$\bone$ 
and that $W\bone = \bone$,}
the detailed formulae become more complicated, but the basic message remains the same; \adb{the details are given in 
Appendix~\ref{4.3-extras}}.  The matrix~$H_\r$ can be represented in the form
\eq\label{H-rho-general}
     H_\r \Eq H_* + \sum_{l=2}^m x_\r\ul (x_\r\ul)\TT,
\en
where $H_* := n^{-1}\bone\bone\TT$ is as before, and $x_\r\ut,\ldots,x_\r\um$ are {\it orthonormal\/} vectors 
orthogonal to~$\bone$ that, with~$\bone$, span the space generated by the columns of $K_\r X$.
The expression for $\t_0^2$ derived in \Ref{tau-0-expanded}, in addition to the quantity
\eq\label{v-rho-def} 
    v(\r_0) \Def \Tr\{(Z_{\r_0}Z_{\r_0}\TT + Z_{\r_0}^2)(I - H_{\r_0})\},
\en 
also involves the traces of the matrices $Z_{\r_0}^2 H_{\r_0}$, 
$Z_{\r_0}\TT Z_{\r_0} H_{\r_0}$, $(Z_{\r_0} H_{\r_0})^2$ and $Z_{\r_0}\TT H_{\r_0} Z_{\r_0} H_{\r_0}$.
\adb{Similarly, from~\Ref{U-dash-expec}, the quantity~$\D_0$ differs from~$\half \t_0^2$ by a quantity that
is also expressed in terms of these traces.}
It is shown in Appendix~\ref{4.3-extras} that if,
for instance, it is assumed that each of the vectors~$x_\r\ul$,
$2\le l\le m$, has components bounded in modulus by~$cn^{-1/2}$, for some fixed~$c \ge 1$ --- 
as is the case for $x\ui = n^{-1/2}\bone$ --- then
\ignore{
, for any $1 \le l,l' \le m$, 
\eqs
    |(x_\r\ul)\TT Z_\r x_\r\uld| 
      &\le& \frac{c^2}{1 - |\r|}\,,
\quad\mbox{and}\quad
    |(x_\r\ul)\TT Z_\r^2 x_\r\uld| 
   \Le \frac{c^2}{(1 - |\r|)^2}\,.
\ens
It then follows that 
}
each of the traces listed above is bounded in modulus by $m^2 c^4/(1-|\r|)^2$.
This gives rise to the following

\msk\nin\adbb{{\bf Conclusion~3:}}
If, \red{in the setting of this section,} $m$ \adb{and~$c$ are also fixed for all~$n$, and~$|\r_0|$ is bounded 
away from~$1$,} then the formulae for $\t_0^2$ and~$\D_0$ given in \Ref{tau0-def} and~\Ref{U-dash-expec}
differ from $4v(\r_0)$ and~$2v(\r_0)$ by quantities that are uniformly bounded in~$n$ \adb{as $n \to \infty$}, 
\adbb{where~$v(\r_0)$ is as defined in~\Ref{v-rho-def}.  The variability
in the MLE of~$\r$ is then on the scale $1/\max\{1,\sqrt{v(\r_0)}\}$.}  

\msk\nin
\adbb{Again, this scale of variation} may be much larger than~$\g_0$
if $v(\r_0)$ is much smaller than $\Tr\{Z_{\r_0}Z_{\r_0}\TT + Z_{\r_0}^2\}$.

However, the principal difference is in the quantity~$\m_0$ that is responsible for
the bias in estimating~$\r$, which, \adbb{from \Ref{mu0-def}, \Ref{rhohat-error}, \Ref{nU-dash-approx} and~\Ref{trace-MH}}, 
\adb{and because 
$(I - \r W)^{-1}\bone = (1-\r)^{-1}\bone$,} can be expressed in the following

\msk\nin\adbb{{\bf Bias formula:}  assuming that $W\bone = \bone$, and that $x_\r\ut,\ldots,x_\r\um$ are the remaining
vectors in the expression~\Ref{H-rho-general} for~$H_\r$, the bias in the MLE is approximately given by} \red{$-\m_0/\D_0$,
where}
\[
     \m_0 \adbb{\Eq 2\Tr\{Z_{\r_0}H_{0}\}} \Eq 2(1-\r_0)^{-1} + \sum_{l=2}^m (x_{\r_0}\ul)\TT Z_{\r_0} x_{\r_0}\ul.
\]

\msk\nin
In the simulations in \cite{mizruchi2008effect}, the remaining columns of~$X$ were realized from independent 
\adb{standard} normally
distributed random variables, and therefore had no particular relationship with $A$ or~$W$, \adb{or with each other.  
\ignore{
Now, for a vector~$\bnu'$ 
of independent \adb{centred} normally distributed random variables, each with variance~$1/\sqrt n$, it follows 
from \Ref{quad-form-mean} and~\Ref{quad-form-var} \adb{(with $C := \half(Z_{\r_0} + Z_{\r_0}\TT)$)} that
\eqs
    \ex\{(\bnu')\TT Z_{\r_0} \bnu'\} &=& n^{-1}\Tr\{Z_{r_0}\} \Eq n^{-1}\sum_{l\ge1}\r_0^{l} \Tr\{W^{l+1}\};\\
    \var\{(\bnu')\TT Z_{\r_0} \bnu'\} &=& n^{-2}\half\Tr\{(Z_{\r_0} + Z_{\r_0}\TT)^2\} \Eq 
           n^{-2}\Tr\{Z_{\r_0}^2 + Z_{\r_0}Z_{\r_0}\TT\}.
\ens
}
In such circumstances, they make no appreciable contribution to the bias.  To see this, consider, for
simplicity, what happens if $\r_0 = 0$.
\ignore{
For $2 \le l\le m$, the column $K_\r x\ul$ is a realization from a multivariate normal distribution 
$\MVN(0,K_\r K_\r\TT)$, and the realizations for different~$l$ are drawn independently.  Now, if $W\bone = \bone$, 
it follows as for~\Ref{W-trace-powers-1} that 
\[
      n \Le \Tr\{K_\r K_\r\TT\} \ =:\ nt_r^2 \Eq n + \r^2\Tr\{WW\TT\} \Le 2n,
\]
and that
\[
       \Tr\{K_\r K_\r\TT K_\r K_\r\TT\}  \Le n(1 + 6\r^2 + 4|\r|^3 + \r^4).
\]
Hence $|K_\r x\ul|^2 = nt_\r^2 + O_P(n^{1/2})$ for each $2 \le l \le m$, and $(K_\r x\uld)\TT K_\r x\ul = O_P(n^{1/2})$
for each $2 \le l' < l \le m$.  Thus the
}
Then, for $2\le l \le m$, the {\it orthonormal\/} vectors~$x_{0}\ul$ can be thought of roughly as being columns of~$X$,
divided by $|x\ul| \approx \sqrt n$ so as to have norm~$1$;
`roughly'  involves neglecting the small adjustments that are needed, so as to be made orthogonal to~$\bone$ and 
to each other.
As a result, using  \Ref{quad-form-mean} and~\Ref{quad-form-var} \adb{(with $C := \half(Z_{\r_0} + Z_{\r_0}\TT)$)}, 
it follows that typical values of $(x_{0}\ul)\TT Z_{0} x_{0}\ul = (x_{0}\ul)\TT W x_{0}\ul$ are roughly centred 
around $n^{-1}\Tr\, W = 0$, with standard deviation $n^{-1}\bigl(\Tr\{\half(W + W\TT)^2\}\bigr)^{1/2}$.
\ignore{
vectors~$(t_{\r_0} \sqrt n)^{-1}K_{\r_0} x\ul$, $2\le l\le m$, which are independent, and have distribution
$n^{-1/2}\MVN(0,t_{\r_0}^{-2}K_{\r_0} K_{\r_0}\TT)$; `roughly'  involves
neglecting the small adjustments that are needed, so as to be orthogonal to~$\bone$ and to each other.
As a result, using  \Ref{quad-form-mean} and~\Ref{quad-form-var} \adb{(with $C := \half(Z_{\r_0} + Z_{\r_0}\TT)$)}, 
it follows that typical values of $(x_{\r_0}\ul)\TT Z_{\r_0} x_{\r_0}\ul$ are roughly centred 
around
$$
    n^{-1}t_{\r_0}^{-2}\Tr\{Z_{r_0}K_{\r_0} K_{\r_0}\TT\} \Eq n^{-1}t_{\r_0}^{-2}\Tr\{W K_{\r_0}\TT\}
                                                          \Eq -\r_0 n^{-1}t_{\r_0}^{-2}\Tr\{W^2\} 
$$ 
from~\Ref{Z-rho-def}, withvariance
standard deviation at most $n^{-1}\bigl\{\half\Tr\{(WK_{\r_0}\TT + K_{\r_0}W\TT)^2\}\bigr\}^{1/2}$.
}
Hence they are} individually negligible when compared to \adb{$\t_0 \asymp \{\Tr\{WW\TT\}\}^{1/2}$, and
thus have no appreciable effect on the bias.  \red{The same is true for their sum, if~$m$ is not large.}
A similar, \red{but more complicated,} argument leads to the same conclusion for $\r_0 \neq 0$.}  This phenomenon is
confirmed in the simulations of \cite{mizruchi2008effect}.  However, if the underlying graph encodes some community
structure, with more edges between members of the same community than outside, and if the regression is designed
to investigate community related effects, then the resulting contributions to~$\m_0$ can be as important as that from
fitting an overall mean, and can have a corresponding effect on the bias.
\adbb{Model~$M2$ of Section~\ref{simul} is of this form.} 

\subsection{Ill conditioned graphs}\label{Ill-conditioned}
\adb{
If the weight matrix associated with the underlying network is ill conditioned, there are consequences for
estimation, both of~$\r$ and of~$\b$.  Taking~$\b$ first,
the distribution of~$\hat\b_{\r_0}$ from~\Ref{NDM-MLE} on the model~$M_{\r_0}$ has covariance matrix 
$\s^2 (X\TT\VV_{\r_0} X)^{-1}$,
as in~\Ref{NDM-MLE-2}, where $\VV_\r := K_\r \TT K_\r$.  It is shown in Appendix~\ref{off-the-model} that
the variance of the linear combination~$a\TT \b$, for any given $m$-dimensional unit vector~$a$, is bounded
above by 
\[
    \frac{\s^2}{\lmin(X\TT \VV_{\r_0} X)} \Le \frac{\s^2}{\lmin(X\TT X) \lmin(\VV_{\r_0})},
\]
\red{where~$\lmin(M)$ and~$\lmax(M)$ denote the smallest and largest eigenvalues of a symmetric matrix~$M$.}
Thus, for fixed~$\s^2$, provided that $\lmin(\VV_{\r_0})$ is not close to zero, 
linear combinations of~$\b$ can be accurately estimated if $\lmin(X\TT X)$ is large, which is the usual
condition for ordinary least squares.  And, even if $\lmin(\VV_{\r_0})$ is small, it may still be the
case that $\lmin(X\TT \VV_{\r_0} X)$ is comparable to  $\lmin(X\TT X)$, if the space spanned by the
columns of~$X$ does not contain vectors that are substantially shrunk by~$\VV_{\r_0}^{1/2}$.
}

\adb{
If~$\r_0$ is not accurately estimated, it will typically be the case that $\hat\b_{\r}$ is used to estimate~$\b$,
for some value $\r \neq \r_0$.  In such cases, it is
shown in Appendix~\ref{off-the-model} that linear combinations of~$\b$ can be estimated using~$\hat\b_{\r}$
to an accuracy comparable to that obtained using~$\hat\b_{\r_0}$,
provided that $\lmax(\VV(\r,\r_0))$ is not large, where
 $\VV(\r,\r_0) := \VV_{\r}^{1/2} \VV_{\r_0}^{-1} \VV_{\r}^{1/2}$.  However, if $\lmax(\VV(\r,\r'))$ can be large
for values of $\r$ and~$\r'$ that are plausible, estimation of~$\b$ may be rather less accurate. In particular,
\red{it is shown in Appendix~\ref{off-the-model} that} a large value of~$\lmax(W\TT W)$ 
indicates that problems with the estimation of~$\b$ may also arise;  \red{see Example~1 in Section~\ref{Ill-conditioned}.}
}

\adb{For bias and variability in the estimation of~$\r$, it has already been shown that problems can
easily arise in dense networks, without~$W$ having any unusual structure.  However, these problems can}
be accentuated by particular networks and design matrices. Again assuming that $W\bone = \bone$,  
this is because the trace $\adb{v(\r_0) = }\Tr\{(Z_{\r_0}Z_{\r_0}\TT + Z_{\r_0}^2)(I - H_{\r_0})\}$ in~\Ref{v-rho-def}
can be very much smaller than $\Tr\{Z_{\r_0}Z_{\r_0}\TT + Z_{\r_0}^2\} = 1/\g_{\r_0}^2$.
If $\t_0 \ge 1$, so that $\ssc_0 := \t_0/\D_0 \asymp 1/(\t_0 + 1)$, this leads to a value of $\ssc_0$ 
that is much larger than the theoretical lower bound~$\g_{\r_0}$. 
Taking $\r_0=0$ for simplicity, the comparison is between $\Tr\{(WW\TT + W^2)(I-H_0)\}$ and $\Tr\{WW\TT\}$.

\msk\nin{\bf Example~1.}
As a first example, 
take that of the $(n-1)$-star, for which $A_{1j} = A_{j1} = 1$, $2 \le j\le n$, and
$A_{ij}=0$ otherwise,  and take $W := D^{-1}A$.
Then the matrix~$W^2$ has $(W^2)_{11} = 1$, and all other elements less than or equal to~$1/(n-1)$.
If \adb{$X$ consists only of the single column~$\bone$,} then $H_0 := H_* := n^{-1}\bone\bone\TT$,  $WH_0 = H_*$, 
and the matrix $W^2 H_0 = H_*$ has all its elements 
equal to~$1/n$.  Thus the traces of~$W^2$ and~$W^2H_0$
are bounded for all~$n$; direct calculation gives $\Tr\{W^2\} = 2$ and $\Tr\{W^2 H_0\} = 1$.  The matrix~$WW\TT$ 
consists of an $(n-1)\times(n-1)$ block of $1$'s, for $2 \le i,j \le n$, together with $1/(n-1)$ at position $(1,1)$,
and the remaining elements zero; the matrix $WW\TT H$ has a corresponding $(n-1)\times(n-1)$ block, now with values $1 - 1/n$,
all elements in the first row have the value $1/n(n-1)$, and the remaining elements in the first column 
have value $1-1/n$.
Thus the diagonal elements of~$WW\TT(I-H_0)$ all have values equal to~$1/n$, and so its trace is exactly~$1$,
whereas $\Tr\{WW\TT\} = n-1+1/(n-1)\ge n-1$.  Hence the lower bound~$\g_0$ at $\r=0$ is less than~$1/(n-1)$, but, with
the simple structure corresponding to estimating the overall mean in the model, \adbb{using the formula~\Ref{with-structure}
from Section~\ref{Estn-with-structure},}
\[
     \t_0^2 \Eq 4\Tr\{(WW\TT + W^2)(I - H_*)\} \Eq 4(1 + 2 - 1) \Eq 8
\]
is comparable to~$1$, implying that the MLE does not give accurate estimation of~$\r$.

\adb{
For estimating~$\b$, note that, from \Ref{NDM-MLE} and~\Ref{NDM-MLE-2}, the variance of~$\hat\b_\r$ 
on the model~$M_{\r}$ is given by
$\s^2(\bone\TT K_\r\TT K_\r \bone)^{-1} = \s^2/\{n(1-\r)^2\}$, because $K_\r \bone = (1-\r)\bone$.
Hence, if~$\r$ is known, $\b$ can be estimated accurately by~$\hat\b_{\r}$ if~$n$ is large and~$\r$ is not too close
to~$1$.  However, in this example, if $e\ui$ denotes the vector \red{$(1,0,\ldots,0)\TT$,} then 
$\lmax(W\TT W) \ge (e\ui)\TT W\TT W e\ui = n-1$ is large, suggesting that using the wrong value of~$\r$ in estimating~$\b$
may cause problems.  Applying \Ref{off-model-3} from Appendix~\ref{off-the-model} below, and taking~$\r_0=0$ to simplify
the calculations, \red{since then $S_{\r_0} = I$,} we obtain the explicit expression
\eqs
    \var_{0} \hat\b_\r &=& \frac{\s^2}{n^2(1-\r)^2} \bone\TT K_\r K_\r\TT \bone 
                       \Eq  \frac{\s^2}{n^2(1-\r)^2} \{n(1-2\r) + \r^2 \bone\TT W W\TT \bone\} \\
                       &=& \frac{\s^2}{n^2(1-\r)^2} \Bigl\{n(1-2\r) + \r^2\Bigl((n-1)^2 + \frac1{n-1}\Bigr)\Bigr\}.
\ens
Hence the variance of~$\hat\b_\r$ grows like~$\s^2\{\r/(1-\r)\}^2$ away from $\r_0 = 0$, and is thus substantial
when the wrong value of~$\r$ is used.  As a result, since~$\r$ cannot be reliably estimated, nor can~$\b$.
}

\msk\nin\adbb{{\bf Conclusion~4:}}
This example shows that having a large number of vertices of low degree --- in this case, $n-1$ vertices of degree~$1$ 
--- is {\it not\/} of itself enough to imply accurate estimation of $\r$ \adbb{or~$\b$}, even though it ensures 
that $\Tr\{WW\TT\}$ is large, and hence that~$\g_0$ is small.  

\msk
Adding further columns to~$X$ results in a projection matrix
$H = H_* + \tH$, with $H_* \tH = 0$, $\tH\TT = \tH$ and $\tH^2 = \tH$, replacing the trace $\Tr\{WW\TT(I-H_*)\}$ by 
$\Tr\{WW\TT(I-H)\}$; the latter is smaller, because 
\eqa
    \lefteqn{\Tr\{WW\TT(I-H_*)\} - \Tr\{WW\TT(I-H)\}} \non\\
    &&\Eq \Tr\{WW\TT(I-H_*)\} - \Tr\{(I-H_* - \tH)WW\TT(I-H_* - \tH)\} \non\\
    &&\Eq \Tr\{WW\TT(I-H_*)\} \non\\
    &&\qquad\qquad\mbox{} - \Tr\{(I-H_*)WW\TT(I-H_*) - 2\tH WW\TT(I-H_*) + \tH WW\TT \tH\} \non\\
    &&\adbb{\Eq \Tr\{2(I-H_*)\tH WW\TT - \tH WW\TT\} \Eq \Tr\{\tH WW\TT\}} \non\\
    &&\Eq \Tr\{(\tH W)(\tH W)\TT\} \ \ge\ 0.  \label{smaller-trace}
\ena
The value of $\Tr\{W^2 H\}$ is not materially changed, because, apart from~$W_{11} = 1$, the elements of~$W^2$ are 
uniformly of order~$O(1/n)$. \adbb{Hence, even with the more realistic model structure, the value of~$v(\r_0)$
is still comparable to~$1$ when $\r_0 = 0$, and, in view of Conclusion~3, the variability of the MLE of~$\r$ is 
substantial.}

\msk\nin{\bf Example~2.}
The same sort of discussion can be carried out, if the network above is generalized to a small number of stars,
each of which has
high degree, the remaining vertices being leaves, and the centres of the stars being connected among themselves.
The elements of~$W^2$, apart from the diagonal elements corresponding to the centres of the stars, 
are then uniformly small, as in~\Ref{W-traces}, typically resulting in only moderate values
of $\Tr\{W^2\}$ and $\Tr\{W^2H\}$.  The block structure of~$WW\TT$ in the previous example is replaced by a collection of
blocks of~$1$'s, one for the leaves of each star. As a result, by arguments similar to those in the example
above, the trace $\Tr\{WW\TT(I-H)\}$ becomes small if the design matrix~$X$ contains  columns identifying the blocks; 
the elements of the column identifying a block should take the value~$1$ at positions within the
corresponding block, and zero within the other blocks, the remaining elements being chosen to ensure that
the sum of the columns is~$\bone$, so that the overall mean belongs to the design. \adbb{In consequence,
the value of~$v(0)$ is comparable to~$1$, and, in view of Conclusion~3, the variability of the MLE of~$\r$ is 
substantial in this case as well.}

\bsk\nin{\bf Example~3.}
A further, more extreme example, is the one considered
by \cite{smith2009estimation}, in which $A = \bone\bone\TT - I$ is the adjacency matrix of the {\it complete\/} 
graph on~$n \ge 2$ vertices. Then, for $W = D^{-1}A$, we have 
\adb{
\eq\label{EX3-0}
   W \Eq W\TT \Eq \frac1{n-1}\{\bone\bone\TT - I\} \Eq H_* - \frac1{n-1}\{I - H_*\},
\en
where, as before, $H_* = n^{-1}\bone\bone\TT$, implying, since $H_*$ and $I-H_*$ are idempotent and $H_*(I-H_*)=0$,  that
\eq\label{Ex3-1}
    W^r \Eq H_* + \Bigl(\frac{-1}{n-1}\Bigr)^r \{I - H_*\} \quad\mbox{and}\quad 
              Z_\r \Eq  \frac{H_*}{1-\r} - \frac{I - H_*}{n-1+\r} \red{\Eq Z_\r\TT}.
\en
Hence,} in particular,
\ignore{
$W = W\TT = (n-1)^{-1}\{\bone\bone\TT - I\}$.
$W$ has eigenvector~$\bone$ with eigenvalue~$1$, as well as the $(n-1)$-dimensional eigenspace orthogonal to~$\bone$,
all of whose vectors have eigenvalue~$-1/(n-1)$, giving $\Tr\{W^r\} = 1 + (-1)^r/(n-1)^{r-1}$, $r \ge 1$, and, in particular,}
$$
     \red{\Tr\{Z_\r Z_\r\TT + Z_\r^2\} = 2\Bigl\{\frac1{(1-\r)^2} + \frac{n-1}{(n-1+\r)^2} \Bigr\},} 
$$
and so, in view of~\Ref{CR}, estimation of~$\r$ cannot be accurate \red{if~$\r$ is not close to~$1$.}  
\adb{However, if the model consists only of an intercept term, so that $X = \bone$, 
\red{and recalling that $\bone\TT W = \bone\TT$,} the estimate of~$\b$ 
from~\Ref{NDM-MLE} is given by
\[
      \hat\b_\r(y) \Eq (\bone\TT (I - \r W)^2 \bone)^{-1} \bone\TT (I - \r W)^2 y
                   \Eq n^{-1}\sum_{j=1}^n y_j,
\]
which is the same for all values of~$\r$, and its variance, obtained using the general formula in~\Ref{off-model-3} below, 
reduces to $1/\{n(1-\r_0)^2\}$ on~$M_{\r_0}$, and becomes small as~$n$ increases.  Thus, although~$\r$ may not
be accurately estimable, there is no problem in estimating~$\b$.
}

\adb{Returning to the estimation of~$\r$, still with $X = \bone$, note also that, using~\Ref{Ex3-1},
\[
   \Tr\{(I-H_0)WW\TT\} \Eq \Tr\{(I-H_*)W^2\} \Eq \frac1{(n-1)^2}\Tr\{I - H_*\} \Eq \frac1{n-1}.
\]
Thus}
$\Tr\{(I-H_0)WW\TT\}$ is an order of magnitude smaller than~$\Tr\{WW\TT\}$.
This in turn, from~\Ref{with-structure}, implies that the variability in the MLE of~$\rho$ is on the scale
of~$\t_0 \ll 1$,  in apparent contradiction to the lower bound in~\Ref{CR};  this can only be the case if
the {\it bias\/} in the MLE is extremely large.
We now investigate what happens in more detail.

First, allowing a more general design matrix~$X$ \adbb{with~$m$ columns,} observe that
\[
      WX \Eq (n-1)^{-1}\{\bone\bone\TT X - X\},
\]
so that the $i$-th column of~$WX$ is a linear combination of~$\bone$ and the $i$-th column of~$X$. If 
 $X$ contains the column~$\bone$,
this implies that the columns of~$K_\r X = X - \r WX$ span the same space as the columns of~$X$, and hence
that the orthogonal projection~$H_\r$ is the same for all~$\r$: $H_\r = H_0 = H(X)$, \adb{and $H_0\bone = \bone$.
As a result, it follows that
\eq\label{Ex3-2}
    H_0 H_* \Eq H_* H_0 \Eq H_*  \quad\mbox{and}\quad (I - H_0)H_* \Eq H_*(I - H_0) \Eq 0,
\en
\ignore{
\eq\label{ADB-I-H-W}
       (I-H_\r)W \adbb{\Eq W - H_0 W} \Eq W - (n-1)^{-1}\{\bone\bone\TT - H_0\} \Eq -(n-1)^{-1}(I - H_0);
\en
applying this repeatedly in turn gives
\eqs
    \lefteqn{\Tr\{(I - H_{0})Z_{\r_0}Z_{\r_0}\TT\} \Eq \adbb{\Tr\{(I-H_0)W(I-\r_0 W)^{-1}W(I-\r_0 W)^{-1}\}} } \\ 
     &&\Eq \adbb{\sum_{r\ge0}(r+1)\r_0^r \Tr\{(I-H_0)W^{r+2}\} 
                  \Eq \frac{\Tr\{I-H_0\}}{(n-1)^2} \sum_{r\ge0}(r+1)\Bigl(\frac{-\r_0}{n-1}\Bigr)^r} \\
     &&\Eq \frac{n-m}{(n-1)^2}\Bigl(1 - \frac{\r_0}{n-1}\Bigr)^{-2}  \ \approx\ \frac1{n}\,,
\ens
}
so that, from \Ref{tau0-def}, \Ref{Ex3-1} and~\Ref{Ex3-2},
\[
   \t_0^2 \Eq 8\Tr\{Z_{\r_0}(I - H_0) Z_{\r_0}(I - H_0)\} \Eq \frac{8}{(n-1+\r_0)^2}\Tr\{I - H_0\} 
              \Eq \frac{8(n-m)}{(n-1+\r_0)^2} \ \approx\ \frac8n\,.
\]
Similar calculations show that, using~\Ref{Ex3-1} and~\Ref{Ex3-2} in~\Ref{U-dash-expec}, there is great simplification,
yielding 
\eqs
     \D_0 &=&  \half\t_0^2 + 2\Tr\{- Z_{\r_0} H_{0}(Z_{\r_0} + Z_{\r_0}\TT)(I - H_{0}) + H_0 Z_{\r_0}^2\} \\
          &=&  \half\t_0^2 + 2\Bigl\{\frac1{(1-\r_0)^2} + \frac{m-1}{(n-1+\r_0)^2} \Bigr\}
                 \adbp{\Eq \frac2{(1-\r_0)^2} + O(n^{-1})},
\ens
since $\Tr\{H_0 - H_*\} = m-1$; and then
\eqs 
    \m_0 &=& 2\Tr\{H_0 Z_{\r_0}\} \Eq \frac2{1-\r_0} \red{- \frac{m-1}{n-1+\r}}.
\ens
}
Hence it follows that \adbp{$\m_0/\D_0$ is not small, and that the bias in the MLE is thus substantial.}  

As it happens, \adb{using \Ref{EX3-0} and~\Ref{Ex3-2}},
\[
      (I - H_\r)K_\r y \Eq (I - H_0)(I - \r W) y \Eq \Bigl\{1 + \frac{\r}{n-1}\Bigr\}(I - H_0) y,
\]
and, from~\Ref{NDM-MLE}, and because $(I-H_\r) = (I-H_\r)^2 = (I - H_0)^2$, this implies that
\[
   \shat^2_\r(y) \Eq n^{-1} (K_\r y)\TT (I - H_\r) K_\r y \Eq  \Bigl\{1 + \frac{\r}{n-1}\Bigr\}^2 n^{-1} y\TT(I - H_0) y.
\]
Thus $\shat^2_\r(y)$ varies very little as~$\r$ ranges over the interval $(-1,1)$.  What is worse,
$\log(\shat^2_\r(y))$ is the sum of a function of~$\r$ alone and a function of the data alone, as is therefore 
 the function~$F(\r;y)$ in~\Ref{infer-disturbance-rho-hat} also.
Thus, in this setting, the maximum likelihood estimator, choosing~$\r \in (-1,1)$ to minimize~$F(\r;y)$, represents a 
choice on which the data have no influence, so that its variability is actually zero, and the same
fixed choice of~$\rho$ would always be made.

\msk\nin\adbb{{\bf Conclusion~5:}}
For the complete graph, the likelihood gives no information
at all about the value of~$\r$.

\msk
\adbp{For this choice of~$W$, $K_\r$ has eigenvalues $1-\r$, together with $1 + \r/(n-1)$ repeated $n-1$ times,
and so the function}
\eqs
     F(\r;y) &=& \log\{n^{-1} y\TT(I - H_0) y\} + 2\log\{1 + \r/(n-1)\} \\
                &&\qquad\mbox{}                - \tfrac2n\,\bigl\{\log(1-\r) +(n-1)\log\{1 + \r/(n-1)\}\bigr\} \\
          &=& \log\{n^{-1} y\TT(I - H_0) y\} + \tfrac2n\,\bigl\{\log\{1 + \r/(n-1)\} - \log(1-\r)\bigr\}
\ens
\red{defined in~\Ref{infer-disturbance-rho-hat}}
\adbb{has strictly positive derivative throughout the (open) interval $(-\infty,1)$. 
Hence setting its derivative equal to zero to find a minimum
would also fail here.}

\section{Estimation using a quadratic form}\label{Quadratic-form-estimation}
\setcounter{equation}{0}

\adbb{The conclusions of Section~\ref{MLE-section} indicate that, in dense graphs, the MLE of the 
correlation parameter~$\r$ 
in the network disturbance model may quite generally exhibit significant bias and variability.  
The ill--conditioned networks of Section~\ref{Ill-conditioned} indicate that these effects can be
dramatic.  Because of the Cram\'er--Rao bound given in~\Ref{CR}, it is not possible to avoid the
variability, whatever estimator of~$\r$ is used.  However, estimators with smaller bias can be envisaged;
see \cite{bao}, \cite{yubaiding} and \cite{yang} for progress in the context of correcting the bias in 
the MLE for network effects models, when the variability is small.
In this section, we move away from the MLE, and instead propose a simple and intuitive estimator of~$\rho$,
which suffers much less from bias, even in circumstances
in which all estimators have substantial variability.}

\subsection{A simple estimator of~$\r$}\label{new-estimator}
An important element in $n\frac{dF(\r;y)}{d\r}$, the quantity that is set equal to zero for the MLE, 
is the quadratic form
\[
    (K_\r y)\TT (I-H_\r)(Z_\r + Z_\r\TT)(I-H_\r) K_\r y
\]
in~$y$, \red{as can be seen in~\Ref{D-F}}.  Here, we replace this quadratic form with
a simpler one.  Since, under the model~$M_{\r_0}$ from~\Ref{M-rho-def}, the 
innovations~$\bnu_i = (K_{\r_0} \beps)_i$, $1\le i\le n$, are independent, 
even at neighbouring vertices, 
\adbb{
the expectation of $(K_{\r_0}\beps)_i (K_{\r_0}\beps)_j$ is zero for all $i \ne j$.
On the other hand, 
for $\r \neq \r_0$, the expectation of $(K_\r \beps)_i (K_\r \beps)_j$ under~$M_{\r_0}$ is typically not zero,
if $i$ and~$j$ are neighbours.  For example, if $\r_0=0$, then, for $i \ne j$, 
\eq\label{QuadForm-W-effect}
     \adb{\s^{-2}}\ex_0\{(K_\r \beps)_i (K_\r \beps)_j\} \Eq \{(I - \r W)(I - \r W)\TT\}_{ij} 
                 \Eq -\r(W_{ij} + W_{ji}) + \r^2 (WW\TT)_{ij},
\en
which is not zero whenever $\r \neq 0$ is small enough, if one of $W_{ij}$ and~$W_{ji}$ is non-zero.  
This suggests building sums
\eq\label{T-def}
     T^C_\r(y) \Def (K_\r \hbe_\r(y))\TT C K_\r \hbe_\r(y)       
\en
of products of the form 
$(K_\r \hbe_\r(y))_i (K_\r \hbe_\r(y))_j$, where $\hbe_\r(y) := K_\r^{-1}(I - H_\r)K_\r y$ is the vector of 
residuals after fitting~$M_\r$, as  derived in~\Ref{residuals}; thus
\eq\label{T-def-2}
   T^C_\r(y) \Eq  y\TT K_\r\TT (I - H_\r) C (I - H_\r)K_\r y.
\en  
\adb{
The matrices~$C$, having non-negative coefficients that do not depend on the data or on~$\r$, 
and with all diagonal elements equal to zero, are to be suitably chosen.}
Any such sum~$T^C_\r(y)$, for fixed $C$ and~$y$, should take values close to zero
when~$\r$ is close to~$\r_0$, and this can potentially be exploited to gain information about the value of~$\r$.
As indicated in~\Ref{QuadForm-W-effect}, the elements of~$C$ should be chosen so as to give larger weight to
pairs $(i,j)$ where $W_{ij} + W_{ji}$ is large, and a natural choice is thus to take $C = \half(W + W\TT)$,
or, equivalently, to take $C = W$.  Another intuitive choice, reflecting the structure of the underlying graph,
rather than the model weights, is to take $C=A$, roughly analogous to the lag-one autocovariance in time series
analysis.
}

 Recalling that, on~$M_{\r_0}$, $\by = X\b + K_{\r_0}^{-1}\bnu$, it follows that, \adb{on~$M_{\r_0}$,}
\[
   \hbe_\r(\by) \Eq K_{\r}^{-1}(I-H_\r)K_\r K_{\r_0}^{-1}\bnu,
\]
since $(I-H_\r)K_\r X = 0$, and hence that, on~$M_{\r_0}$,
\eq\label{A-QF}
    T^C_\r(\by) \Eq \bnu\TT Q_C(\r,\r_0) \bnu,
\en
where, much as in~\Ref{Q-def},
\eq\label{QA-def}
    Q_C(\r,\r_0) \Def (K_\r K_{\r_0}^{-1})\TT (I-H_\r) C (I-H_\r) K_\r K_{\r_0}^{-1}.
\en
From \Ref{A-QF} and~\Ref{QA-def}, we have
\eq\label{mu-hat-def}
     \ex_{\r_0} T^C_{\r_0}(\by) \Eq \s^2\Tr\{(I-H_{\r_0}) C\} \Eq -\s^2\Tr\{H_{\r_0}C\} \ =:\ \s^2\hmu^C_{\r_0},
\en
since~$C$ has zero diagonal, and so $\Tr\{C\}=0$.
Equation~\Ref{mu-hat-def} suggests the following `method of moments' procedure for estimating~$\r$.

\bigskip\nin\adbb{{\bf Elementary estimate of}~$\r$:  find $\r = \hrho_C$ so that
\eq\label{estimating-equation}
     T^C_{\r}(y) + \shat_\r^2(y)\Tr\{H_{\r}C\} \Eq 0.
\en
}

\nin Equation~\Ref{estimating-equation} is used to
correct the bias that would result, in view of~\Ref{mu-hat-def}, if, in estimating~$\r$, $T^C_\r(y)$ were itself
set equal to zero. 

\bigskip\nin\adbb{{\bf Procedure}:
\begin{description}
\item[Input:] a network with weight matrix~$W$, a linear model with design matrix~$X$ and observed data~$y$,
               and a non-negative coefficient matrix~$C$ with zero diagonal.
\item[Function:] for any value of~$\r$:
\begin{enumerate}
 \item  Calculate the trace $\Tr\{H_\r C\}$, \adb{using \Ref{H-rho-def} and~\Ref{H-proj-def};}
 \item  Calculate the parameter estimates $\bhat_\r(y)$ and~$\shat^2_\r(y)$ using~\Ref{NDM-MLE};
 \item  Calculate~$\hbe_\r(y)$, using $\bhat_\r(y)$ and~\Ref{residuals};
 \item  Calculate the value of~$T^C_{\r}(y)$ from~\Ref{T-def}.
\end{enumerate}
\item[Solution:]  use a root finding algorithm to find the zero of $T^C_{\r}(y) + \shat_\r^2(y)\Tr\{H_{\r}C\}$,
  which is the estimate~$\hrho_C$ of the correlation~$\r$ using $T^C$.
\end{description}
}

\adb{In this paper, we do not attempt to give detailed conditions under which the procedure
using~\Ref{estimating-equation} gives good estimates $\hrho_C$ of~$\r$.  It is clear from the example of 
the complete graph in Section~\ref{Ill-conditioned} that this
cannot always be so; indeed, for the complete graph, \Ref{estimating-equation} has no solution.  Instead,
we carried out some simulations, illustrating that~$\hrho_C$ can be effective,
within the limitations of what is possible. These are discussed in Section~\ref{simul} below.
In the remainder of this section, we give a rough guide to its precision from a theoretical standpoint,
as well as proposing a permutation procedure for assessing its precision based on the data to hand.}

\msk\nin\adb{{\bf The scale of variation of~$\hrho_C$}}\\
\adbb{
To determine the scale of variation about~$\r_0$ of the resulting estimator~$\hrho_C$, write
\eq\label{U^C-def}
    U^C_\r(y) \Def T^C_{\r}(y) + \shat_\r^2(y)\Tr\{H_{\r}C\}.
\en
As in the discussion of the MLE \adb{around \Ref{rhohat-error} and~\Ref{nU-dash-approx}} in Section~\ref{MLE-section}, 
a plausible measure is given by
$\ssch^C_{\r_0} := \htau^C_{\r_0}/|\hD^C_{\r_0}|$, to be estimated by $\htau^C_{\hrho_C}/|\hD^C_{\hrho_C}| $, where
\[
    (\htau^C_{\r_0})^2 \Def \s^{-4}\var_{\r_0}(U^C_{\r_0}(\by)) \quad\mbox{and}\quad 
                \hD^C_{\r_0} \Def \adb{-}\s^{-2}\ex_{\r_0}\Bigl\{\frac{dU^C_{\r}(\by)}{d\r}\Bigr|_{\r=\r_0}\Bigr\}.
\]
The quantity~$(\htau^C_{\r_0})^2$ is dominated by the trace of a matrix,
\eq\label{tau-hat-def-approx}
   \Tr\{(I-H_{\r_0})(C + C\TT) (I-H_{\r_0})(C + C\TT)\},
\en
which is reminiscient of the definition of~$\t^2_0$ in~\Ref{tau0-def},
and the quantity~$\hD^C_{\r_0}$ has leading term 
\eq\label{Delta-hat-def-approx}
    \Tr\bigl\{(I-H_{\r_0})\adbb{[Z_{\r_0}\TT(I - H_{\r_0}) - Z_{\r_0}H_{\r_0}]} (C + C\TT) \bigr\}.
\en
The derivations, and more detailed expressions, are given in Section~\ref{scale-quadratic-estimator}.
As before, the scale of variation 
\eq\label{var-scale} 
            \ssch^C_{\r_0} \Def \htau^C_{\r_0}/|\hD^C_{\r_0}|,
\en 
determined from \Ref{tau-hat-def-approx} and~\Ref{Delta-hat-def-approx},
need not correspond precisely to the standard deviation of~$\hrho_C$ on~$M_{\r_0}$,
because~$U^C_\r(y)$ is not linear as a function of~$\r$, and because the derivative
of~$U^C_\r(y)$ may fluctuate appreciably about its expectation.
}

\adbb{To get an idea of the magnitude of~$\ssch^C_{\r_0}$, let $C = W$, and
suppose that $W\bone = \bone$ and that $X = \bone$, so that only the overall mean is being fitted. Then,
\adb{as in Section~\ref{Estn-with-structure},}
$K_\r X = (1-\r)\bone$, implying that, for all~$\r$, $H_\r = H_* := n^{-1}\bone\bone\TT = WH_\r$ and that 
$Z_\r H_\r = (1-\r)^{-1}H_*$.  As a result, we easily deduce that
\[
     (\htau^W_{\r_0})^2 \Eq \adb{2}(\Tr\{W^2 + WW\TT - H_* WW\TT\} - 1) + O(n^{-1}),
\]
and that
\[
    \hD^C_{\r_0} \Eq -\Tr\{WZ_{\r_0} + WZ_{\r_0}\TT - H_* WZ_{\r_0}\TT\} + (1-\r)^{-1} 
                + n^{-1}\Tr\{(Z_{\r_0} + Z_{\r_0}\TT)(I - H_*)\},
\]
where we have used the full expressions in \Ref{tau-hat-def} and~\Ref{Delta-hat-def} instead of
\Ref{tau-hat-def-approx} and~\Ref{Delta-hat-def-approx} to obtain the smaller order terms.
In particular, if $\r_0 = 0$, so that $Z_{\r_0} = W$, $\adb{\half}(\htau^W_{0})^2 = -\hD^C_{0} + O(n^{-1})$,
and the magnitude of~$\ssch^C_{0}$ is roughly $2/\htau^W_{0}$;  this remains the case for~$\r_0$
not too far from zero.   As in the discussion in Section~\ref{Cramer-Rao}, this implies that if, for instance, 
$W_{ij} \le c/n$ for all $i,j$, then the scale of variation is at least $1/\sqrt c$, and is appreciable
if~$c$ is not large.  For $W := D^{-1}A$, this corresponds to having a dense underlying graph, as before.
}

The justification for using~\Ref{estimating-equation} to estimate~$\rho$ relies only on the first two joint moments
of the vector~$\bnu$ of errors, and does not assume that the components are independent or normally distributed, though they
should be uncorrelated.  However, the formula~\Ref{tau-hat-def}, used in the theoretical discussion of the scale of
random variation~$\ssch^C_{\r_0}$ of the estimator~$\hrho_C$, is already based on stronger assumptions.  For instance, 
if the errors~$\bnu$ are independent
and have zero mean, but if $\ex\{\s^{-4}\bnu_1^4\} \neq 3$, the expression~\Ref{quad-form-1} for the variance of
a quadratic form \adbb{needs to} be used in place of~\Ref{quad-form-var}. 

\msk\nin\adb{{\bf A permutation procedure for assessing the accuracy of~$\hrho_C$}}\\
\adbb{Instead of assessing the precision of~$\hrho_C$ using~$\ssch^C_{\r_0}$, some form of simulation procedure can be used.
One possibility is the following permutation algorithm.
Start by computing the esimated $\bnu$-residuals
\eq\label{P1}
     \hnu(y) \Def  K_{\hrho_C} \hbe_{\hrho_C}(y) \Eq (I - H_{\hrho_C}) K_{\hrho_C} y.
\en
Since the components of~$\bnu$ are independent and identically
distributed, this suggests that synthetic data could be derived by replacing the $\bnu$-residuals $(\hnu_i(y),\,1\le i\le n)$
by permuted residuals
\eq\label{P2}
    \hnu^\p_i(y) \Def \hnu_{\pi(i)}(y),\qquad 1\le i\le n,      
\en
where~$\pi$ is a randomly chosen permutation of $\{1,2,\ldots,n\}$.  These can be used to determine residuals 
\eq\label{P3}
     \hbe^\pi(y) \Def K_{\hrho_C}^{-1} \hnu^\p
\en
for the original model, and hence to generate synthetic data $y^\pi := X\bhat + \hbe^\pi(y)$.  These can in turn be
analyzed in the same way as the original data, to provide an estimate $\hrho_C^\p$ of~$\hrho_C(y)$.
Repeating this procedure for a large number~$M$ of independently chosen random permutations $\pi_j$, $1\le j\le M$,
yields a selection of estimates~$\hrho_C^{\p_j}$ of~$\hrho_C(y)$, whose empirical variability
should be much like that of the variability of~$\hrho_C$ about \adb{the model value~$\r_0$, when $\r_0 = \hrho_C(y)$.}
}

\subsection{Simulations}\label{simul}
\adb{In this section, we  illustrate the effectiveness of~$\hrho_C$ using simulated data,
showing that it substantially reduces the bias exhibited by the MLE, while its precision
remains constrained by the lower bound~\Ref{CR}.}  In Table~\ref{table1}, data were simulated at the vertices
of a Bernoulli graph~$G(n,p)$, with $n=100$ and with three choices of~$p$.  The data were simulated from the
model~$M_\r$ with six different values of~$\r$. The parameters $\b = (1,0.5,0.4,0.3)$ were held fixed,
and the design matrix~$X$ was chosen to have~$\bone$ as its first column,
with its remaining columns consisting of independent realizations of standard normal random variables. 
For each combination of $p$ and~$\r$, the elements of~$X$ and the network~$G$ were sampled anew; then 100 replicates
were generated,  starting from realizations of independent standard normal random variables~$\bnu$, and
using~\Ref{NDM-errors} to determine the values of~$\e$, \adbp{taking $W_{ij} := D_i^{-1}A_{ij}$ if
$D_i \ge 1$, and $W_{ij} = 0$ if~$i$ is an isolated vertex.}
Each cell in Table~\ref{table1} contains the mean of the corresponding 100 values of~$\hrho_W$, 
\adb{calculated using~\Ref{estimating-equation} with $C = W$,}  together with the standard
error of a single estimate.  It is clear that the standard errors are substantial. The means are close to
the truth, within the limits of variability; for $100$ replicates, the standard error of the mean is one tenth 
of that of a single estimate.  The estimates of~$\beta$ were uniformly good, and are not presented.
In a Bernoulli graph~$G(n,p)$, the degrees of the vertices are reasonably homogeneous, being concentrated
around their mean $(n-1)p$ when~$p \gg 1/n$.
Thus, from~\Ref{W-traces}, calculating as if all the degrees were exactly $(n-1)p$,
the lower bound~\Ref{CR} on the variance of an estimator when $\r=0$ gives~$\g_\r$ approximately equal to~$\sqrt{p/2}$.
The actual values of~$\g_\r$ in the simulations depend on the value of~$\r$ and on the particular realization 
of the Bernoulli random graphs that were sampled. The value of~$\r$ made little difference, but the randomness
in the degrees was still enough, for $n=100$, to reduce the value of~$\g_\r$ by about 10\%; typical values are
quoted in the final line of Table~\ref{table1c}.
They indicate that the observed 
variablility of the quadratic form estimator is not far from the best that could be achieved in this setting.

\begin{table}[ht]
\centering
\begin{subtable}{0.8\textwidth}
  \centering
  \begin{tabular}{c c | c  c  c}\\ 
    \multicolumn{5}{c}{\mbox{}\qquad\qquad\quad$p$} \\[1ex]
      &  & 0.0975 & 0.19 & 0.36 \\ \hline 
      \rule{0pt}{12pt} &  -0.2 & -0.21 (0.20) & -0.24 (0.30) & -0.24 (0.44)     \\   
      &  -0.1 & -0.14 (0.23) & -0.10 (0.35) & -0.11 (0.46) \\  
      $\r$ & 0.0 & 0.01 (0.23) & -0.03 (0.37) & 0.02 (0.47) \\ 
      & 0.1 & 0.08 (0.24)  & 0.09 (0.31) & 0.08 (0.46) \\  
      & 0.2 & 0.21 (0.23) & 0.16 (0.34) & 0.11 (0.50)\\
      & 0.3 & 0.28 (0.21) & 0.23 (0.32) & 0.27 (0.47) \\ \hline 
     \adbb{CRLB} \rule{0pt}{13pt} & 0.0  & (0.20) & (0.29) & (0.39)
  \end{tabular}
  \caption{Means of 100 simulated estimates of $\rho$ (with the standard error of a single estimate) calculated
   using~\Ref{estimating-equation} \adbb{with $C=W$.}  The \adbb{Cram\'er--Rao Lower Bound (CRLB)} line 
   gives the values of the lower bound~$\g_{\r_0}$ derived using~\Ref{CR}.} 
    \label{table1c}
\end{subtable}\\ [2ex]
\begin{subtable}{0.8\textwidth}
   \centering
  \begin{tabular}{c c | c  c  c}\\ 
  \multicolumn{5}{c}{\mbox{}\qquad\qquad\quad$p$} \\[1ex]
   &  & 0.0975 & 0.19 & 0.36 \\ \hline 
   \rule{0pt}{12pt} &  -0.2 & -0.26 (0.21) & -0.28 (0.28) & -0.33 (0.39)      \\   
   &  -0.1 & -0.13 (0.20) & -0.22 (0.34) & -0.28 (0.37)  \\  
   $\r$ & 0.0 & -0.03 (0.21) & -0.10 (0.30) & -0.19 (0.42)  \\ 
    & 0.1 & 0.02 (0.21) & -0.05 (0.30)  & -0.11 (0.41) \\  
   & 0.2 & 0.13 (0.20) & 0.14 (0.31)  & -0.05 (0.42) \\
   & 0.3 & 0.19 (0.20) & 0.16 (0.30)  & 0.02 (0.33)  \\ \hline 
   \adbb{CRLB} \rule{0pt}{13pt} & 0.0  & (0.20) & (0.29) & (0.39)
   \end{tabular}
   \caption{Means of 100 simulated MLE estimates of $\rho$. The standard errors of a single observation are 
    consistent with the values given in the CRLB line.} \label{table1ml}
\end{subtable}
\caption{Comparison between $\hrho_C$ and~$\hrho\ML$ for a network disturbance model based on $G(100,p)$ 
       \adbb{with $W = D^{-1}A$}.}
\label{table1}
\end{table}

For comparison, we give the results in Table~\ref{table1ml} of using the MLE in the same setting.  
The standard deviations were
uniformly close to the values given in the CRLB line in the table above, and were thus a little smaller than
those in the experiments using~\Ref{estimating-equation} \adbb{with $C = W = D^{-1}A$}.
However, the MLE exhibits a systematic bias, which becomes progressively worse
as the density of the graph increases.  Even with edge density $p = 0.0975$, the smallest value investigated,
the means of the estimates of~$\r$ were consistently too small, by amounts comparable to~$0.05$, whereas
the whole range of values of~$\r$ being investigated was only~$0.5$.

\ignore{
\begin{table}[ht]
\begin{center}
\begin{tabular}{c c | c  c  c}\\ 
  \multicolumn{5}{c}{\mbox{}\qquad\qquad\quad$p$} \\[1ex]
 &  & 0.0975 & 0.19 & 0.36 \\ \hline 
\rule{0pt}{12pt} &  -0.2 & -0.26 (0.21) & -0.28 (0.28) & -0.33 (0.39)      \\   
 &  -0.1 & -0.13 (0.20) & -0.22 (0.34) & -0.28 (0.37)  \\  
$\r$ & 0.0 & -0.03 (0.21) & -0.10 (0.30) & -0.19 (0.42)  \\ 
 & 0.1 & 0.02 (0.21) & -0.05 (0.30)  & -0.11 (0.41) \\  
 & 0.2 & 0.13 (0.20) & 0.14 (0.31)  & -0.05 (0.42) \\
 & 0.3 & 0.19 (0.20) & 0.16 (0.30)  & 0.02 (0.33)  \\ \hline 
 \adbb{CRLB} \rule{0pt}{13pt} & 0.0  & (0.20) & (0.29) & (0.39)
\end{tabular}
\caption{Means of 100 simulated MLE estimates of $\rho$, for a network disturbance
model based on $G(100,p)$ \adbb{with $W := D^{-1}A$}. The standard errors of a single observation are consistent with the
values given in the CRLB line.} \label{table1n}
\end{center} 
\end{table}
}

In the second experiment, the graph $G(n,p)$ was replaced by an Erd\H os--R\'enyi mixture model, consisting of
two classes with~$50$ vertices in each;  the within class probabilities were taken to be $2p(1-p)$ and the
between class probabilities to be~$p$, for $p = 0.05,0.1$ and $0.2$.  The value of~$\r$ was kept constant \adbb{at~$0.1$,
and the} parameters~$\b$ were as in the first experiment. \adbb{The results of estimation using~$\hrho_C$
are given in Table~\ref{table2c}.}
In the first row, model~$M1$, the design matrix was also constructed as before, with 
\adb{the first column taken to be~$\bone$, and with the remaining columns} chosen 
to have \adb{independent standard normal}
entries.  For the second row, model~$M2$, the second column of the design matrix
had elements $+1$ for the vertices in the first of the classes, and $-1$ for vertices in the second,
corresponding to an inter-block contrast, which significantly changes the bias correction $\Tr\{H_{\r}C\}$.  In each cell, 
the results of \adbb{1,000} replicates are given.  At this accuracy, a small negative bias in the estimates of~$\r$ 
using~$\hrho_C$ can
be discerned, although it is negligible when compared to the standard errors of the estimates.
The approximate ideal standard errors derived from \Ref{CR} for $\r = 0.1$ are given in the final row.
\adbb{In Table~\ref{table2ml}, the results of using maximum likelihood are given for comparison; here, the bias
is clearly visible.}

\begin{table}[ht]
\centering
\begin{subtable}{0.8\textwidth}
\centering
\begin{tabular}{ c | c  c  c}\\ 
  $p$ & 0.05 & 0.1 & 0.2 \\ \hline 
\rule{0pt}{13pt} $M1$  &  0.102 (0.20) & 0.104 (0.28)  &  0.111 (0.43)   \\   
 $M2$  &  0.094 (0.21)  & 0.082 (0.29) & 0.098 (0.42) \\ \hline
\rule{0pt}{13pt} \adbb{CRLB} & (0.18) & (0.26) & (0.35)
\end{tabular}
\caption{Means of \adbb{1,000} simulated estimates of $\r$ (with the standard error of a single estimate)
using~\Ref{estimating-equation} \adbb{with $C = W$.} 
 The CRLB line gives the values of the lower bound~$\g_{\r_0}$ derived using~\Ref{CR}. } \label{table2c}
\end{subtable}\\[4ex]
\begin{subtable}{0.8\textwidth}
\centering
\begin{tabular}{ c | c  c  c}\\ 
  $p$ & 0.05 & 0.1 & 0.2 \\ \hline 
\rule{0pt}{13pt} $M1$  &  0.06 (0.18) & 0.01 (0.26)  &  -0.09 (0.37)   \\   
 $M2$  &  0.05 (0.19)  & -0.02 (0.27) & -0.14 (0.37) \\ \hline
\rule{0pt}{13pt} CRLB & (0.18) & (0.26) & (0.35)
\end{tabular}
\caption{
\adbb{Means of 1,000 simulated estimates of $\r$ (with the standard error of a single estimate)
using the MLE. The standard errors of a single observation are consistent with the
values given in the CRLB line.} } \label{table2ml}
\end{subtable}
\caption{Comparison between $\hrho_C$ and~$\hrho\ML$ for two network disturbance models based on an Erd\H os--R\'enyi
   mixture graph, each with \adbb{$W := D^{-1}A$ and} with $\r = 0.1$.}
\label{table2}
\end{table}

\adbb{
Table~\ref{table5} gives the results of assessing the variability of~$\hrho_W$ using the permutation algorithm based on 
\Ref{P1}--\Ref{P3}. The experiments were carried out
in the context of the  network disturbance model~$M2$ based on an Erd\H os--R\'enyi mixture graph.
\adbb{In order to make a sensible comparison with the variability
of the estimates~$\hrho_W$ derived from repeatedly simulating data from the model,  
 it was necessary to apply the permutation algorithm to a realization of the model for which the original 
estimate~$\hrho_W$ was
reasonably close to the true value $\r_0 = 0.1$.  Because of the substantial variability in~$\hrho_W$, this was not the
case in most realizations.}  In each experiment, the simulated variability of~$\hrho_W^{\pi}$
agreed well with the variability of~$\hrho_W$ under the model, and the value of~$\ssch^W_{\r_0}$ was also quite similar.
}

\begin{table}[ht]
\begin{center}
\begin{tabular}{ c | c  c  c  c  c}\\ 
  $p$ & $\hrho_W$ & $E(\hrho_W^{\p})$ & $SE(\hrho_W^{\p})$ & $\ssch^W_{\r_0}$ & $SE(\hrho_W)$\\ \hline 
\rule{0pt}{13pt} 
 0.05  &  0.094  & 0.092   &  0.19  &  0.17  &  0.21 \\  
 0.1   &  0.107  & 0.112   &  0.29  &  0.26  &  0.29 \\ 
 0.2   &  0.091  & 0.073   &  0.43  &  0.42  &  0.42  \\ \hline
\rule{0pt}{13pt} 
\end{tabular}
\caption{\adbb{Means $E(\hrho_W^{\p})$ and standard errors $SE(\hrho_W^{\p})$ of 1,000 simulated estimates of $\hrho_W$,
using the permutation procedure, for the model~$M2$ with $\r = 0.1$, together with the scale of 
variation~$\ssch^W_{\r_0}$, calculated using 
\Ref{var-scale}, \Ref{tau-hat-def} and~\Ref{Delta-hat-def}, and the corresponding standard errors $SE(\hrho_W)$ of
$\hrho_W$ from Table~\ref{table2}.} } \label{table5}
\end{center} 
\end{table}

\begin{figure}[ht]
\centering
   \begin{subfigure}[b]{0.4\textwidth}
     \includegraphics[width=\textwidth]{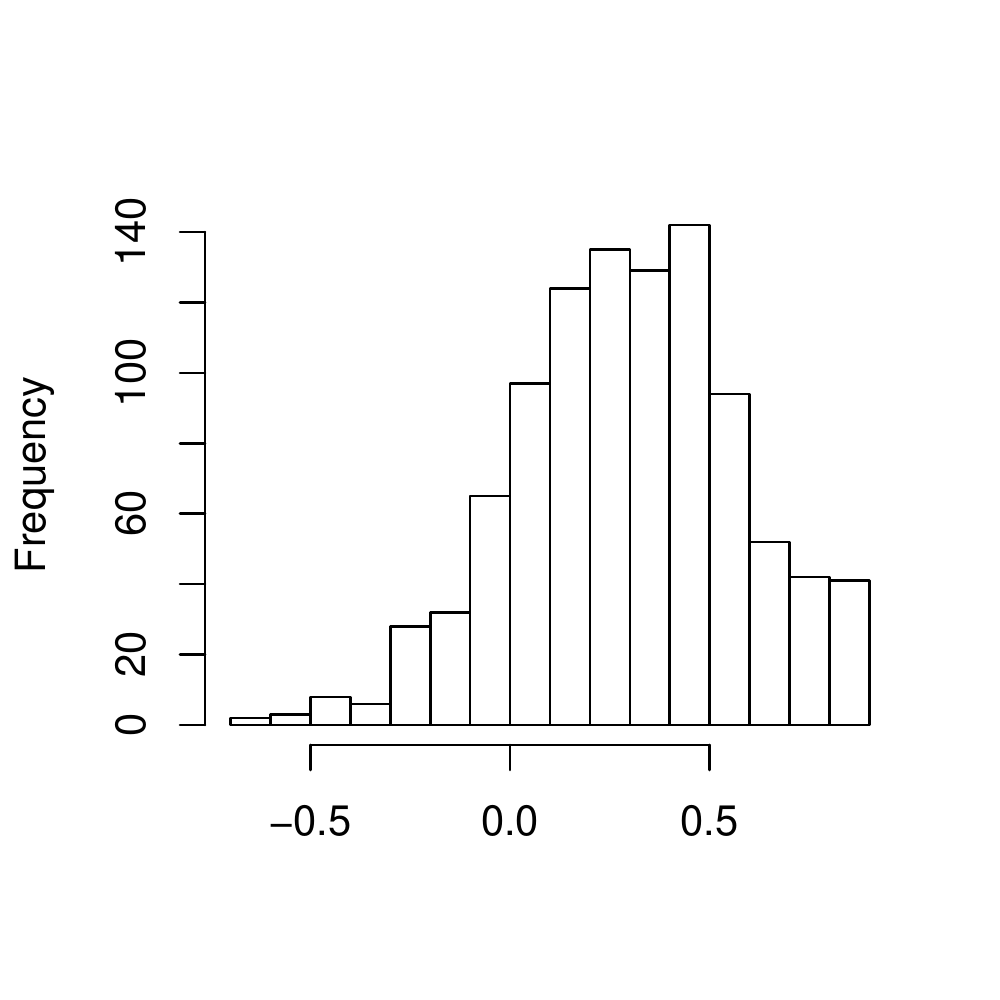}
     \caption{Histogram of $\hrho_W$}
     \label{rhohat-pic}
   \end{subfigure}
   \qquad
   \begin{subfigure}[b]{0.4\textwidth}
     \includegraphics[width=\textwidth]{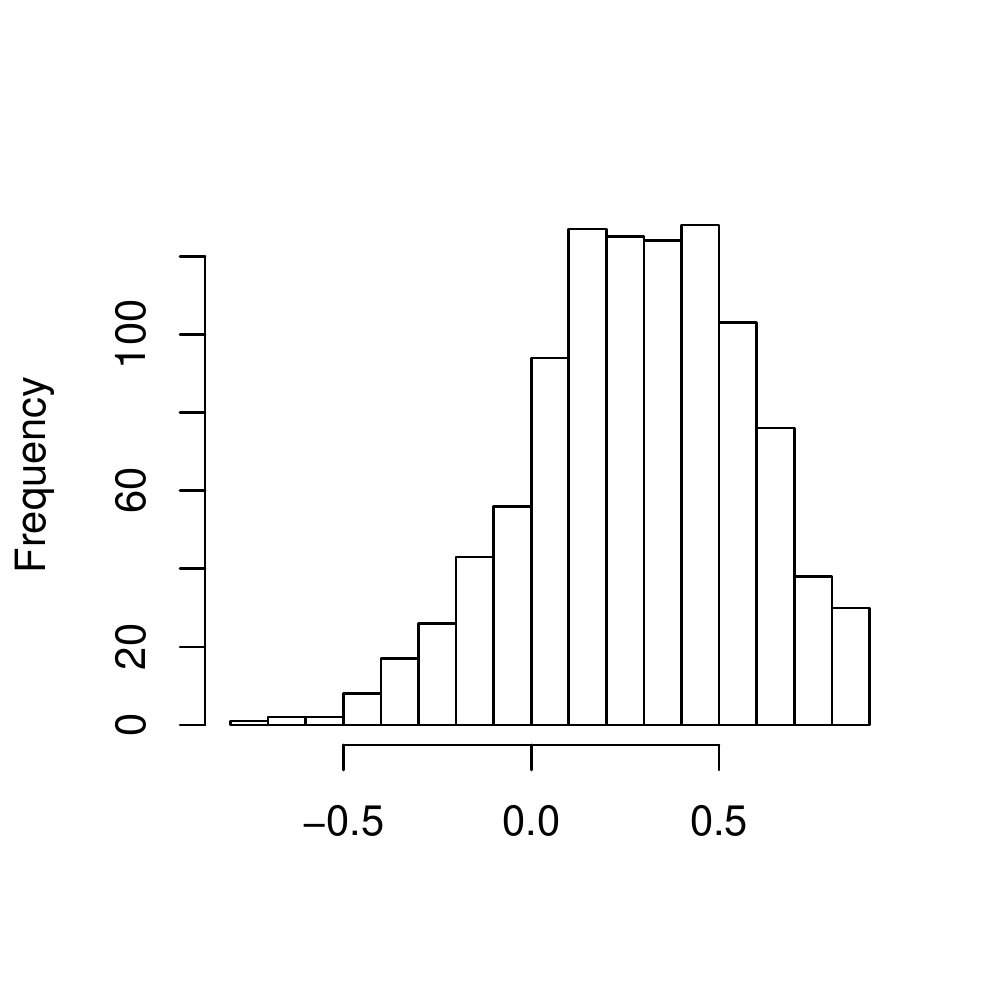}
     \caption{Histogram of $\hrho_W^\p$}
     \label{rhopi-pic}
   \end{subfigure}
\caption{Histogram of $1,000$ values of~$\hrho_W$ for the model~$M2$ with $\rho = 0.3$, together with a histogram
          of $1,000$ values of the permutation estimates~$\hrho_W^\p$ derived from a realization of the same model 
          in which $\hrho_W = 0.302$}
\label{rho-histograms}  
\end{figure}

\adbb{
In Figure~\ref{rho-histograms}, a histogram of $1,000$ values of~$\hrho_W$, simulated from the model~$M2$, 
\adb{{\it but now with\/} $\r = 0.3$,}
is shown next to a histogram of $1,000$ values of~$\hrho_W^\p$, generated using the permutation method described above
from a realization of the model~$M2$ in which $\hrho_W = 0.302$.  The two histograms are very similar, and
neither is substantially skewed, even though the value $0.3$ \adbb{chosen for~$\r$} is close enough to
the upper limit of~$1$ for some compression in the upper tail to be expected.
}

\section{Network effects models}\label{NE-models-sect}
\adbb{
The network effects model~\Ref{NEffects-model} can be written in the form
\eq\label{NEffects-model-1a}
    y \Eq K_\r^{-1}X\b + \e,
\en
where the errors~$\e$ in~$y$ have the same dependence as in the network disturbance model, but the 
{\it model structure\/} explaining~$y$ also changes with~$\r$. 
\adb{  
The corresponding log-likelihood, as in \cite[Equation~(1.7)]{doreian1989network}, is given by
\eq\label{NEM-LL}
   \red{l(\th;y) \Def} -\tfrac n2 \{\log(\s^2) + \log (2\pi)\} - \frac1{2\s^2} (K_\r y-X\b)\TT  (K_\r y-X\b) 
               + \log\det K_\r ,
\en
}
\red{where $\th := (\b,\s^2,\r)$.}   Since the structural element
also contains information about~$\rho$, this is reflected in the estimation procedures.  In particular,
the estimate of~$\b$ using maximum likelihood, for a given value of~$\rho$, is
\[
     \bhat_\r(y) \Def (X\TT X)^{-1}X K_\r y,
\]
as can be deduced from ordinary least squares~\Ref{OLS-MLE} applied to the data $y' := K_\r y$.  As a result, if the
true model is $M_{\r'}$, for $\r' \neq \r$, we have
\eq\label{NEM-beta-mean}
    \ex_{\r'} \bhat_\r(\by) \Eq (X\TT X)^{-1}X K_\r K_{\r'}^{-1}X \b.
\en
This expectation is in general different from~$\b$, in which case the MLE of~$\b$ is likely to be biased if~$\r$ is 
incorrectly estimated.    
Hence, in order to get good estimates of the structural component~$\b$, it is important to have
an accurate estimate of~$\r$.
}

\adbb{
Now, \adb{differentiating~\Ref{NEM-LL} with respect to~$\r$,} the $\r$-component of the score function analogous 
to~\Ref{urho-def}  takes the form
\eqa
   \ts\urho(\th;y) 
      &=& \s^{-2} (K_\r y-X\b)\TT  Z_\r K_\r y - \Tr\{Z_\r\}, \label{urho-def-2}
\ena
where $Z_\r := WK_\r^{-1}$.  Its variance at~$\th_0$ under
the model with $\th = \th_0$ is now given by 
\eq\label{var-score-1}
   V_1 \Def \var_{\th_0}\, \ts\urho(\th_0;\by) 
                \Eq \Tr\{Z_{\r_0}^2 + Z_{\r_0} Z_{\r_0}\TT\} + \s_0^{-2}|Z_{\r_0} X \b_0|^2,
\en
containing an extra term, as compared to~\Ref{var-score}, that involves the parameter~\adbb{$\b_0/\s_0$}.  
This implies that
the lower bound~$1/V_1$ for the variance of an unbiased estimator is smaller than in the analogous network 
disturbance model.
This does not, of itself, mean that estimation of~$\r$ can be made more precise, since there is no guarantee
that the lower bound can be attained; instead, a closer analysis is needed.
}

\adbb{
As an example, suppose that~$X$ consists solely of the column~$\bone$, with parameter~$\b$, and that 
the weight matrix~$W$ satisfies
$W\bone = \bone$.  Then the additional term is just $n(\b_0/(1-\r_0)\s_0)^{2}$, a quantity that grows with~$n$.
On the other hand, the network effects model~\Ref{NEffects-model-1a} translates into
\[
    y \Eq (1-\r)^{-1}\b\bone + \e \Eq \tbe\bone + \e,
\]
which can also be viewed as a network disturbance model, with a new structural parameter~$\tbe = \b/(1-\r)$ in place of~$\b$.  
Thus~$\tbe$ can be accurately estimated, if the network is large enough, as can be deduced from~\Ref{NDM-MLE-2},
but, as before, any estimator of~$\r$ exhibits
substantial variability if the network is dense.   This explains the simulation results observed 
in \cite{mizruchi2008effect}, in which 
the maximum likelihood estimators of both~$\r$ and the coefficient of the constant exhibited  substantial
variability, and were strongly negatively correlated.  In the same experiments, the remaining~$m-1$ columns~$X'$ of~$X$ 
had independent,
normally distributed entries with mean zero and variance~$1$. In this case, with~$\b_1/(1-\r)$ replaced by
a parameter~$\tbe$ and with~$\b' := (\b_2,\ldots,\b_m)^T$, the extra term in the variance of the score function reduces to 
\[
      \s^{-2}\ex\{|Z_{\r} X' \b'|^2\} \Eq (|\b'|/\s)^2 \Tr\{Z_{\r}Z_{\r}\TT\}
\]
(where the expectation is with respect to the random choice of~$X'$), and, \adb{as in~\Ref{W-trace-powers-1.5}, this does not
become large with~$n$ for uniformly dense graphs.}  Hence there are still substantial limitations on the precision of estimation of~$\r$,
for Cram\'er--Rao reasons, and, in consequence, on the estimation of $\b_1 = \tbe(1-\r)$ also.  
}

\adbb{
The same phenomenon makes itself felt if the model structure~$X$ is more complicated, but the columns of~$WX$
all belong to the linear span of the columns of~$X$; $WX = X \G$, for an $m \times m$ matrix~$\G$.
In this case, the model structure~$K_\r^{-1}X\b$ can be expressed as
\eqa
     K_\r^{-1}X\b &=& (I - \r W)^{-1}X\b \Eq X\b + \sum_{l \ge 1} \r^l W^l X\b  \non\\
                  &=& X\b + \sum_{l \ge 1} \r^l X \G^{l}\b \Eq X(I - \r\G)^{-1}\b \ =:\ X\parg,
                       \label{KXb-rephrase}
\ena
where the new parameter $\parg := (I - \r\G)^{-1}\b$.  Once again, this yields a network disturbance
model, in which the parameter~$\parg$ can be accurately estimated if~$n$ is large, whereas, in dense networks,
the estimator of~$\rho$ may exhibit substantial variability.  For the original parameter~$\b$, this would
express itself in variable, but strongly dependent, estimators of the components of~$\b = (I - \r\G)\parg$, which are  
linked through the common value of the estimator of~$\r$.  
}

\adb{
In general, $WX$ can be written in the form
\eq\label{WX-decomposition-1}
     WX \Eq  X\G + E,
\en
where~$E$ is an $n \times m$ matrix, whose columns are orthogonal to those of~$X$,
which measures how far~$WX$ lies from the linear span of the columns of~$X$. 
In practice, except for the model in which only the overall mean is to be estimated,
and $W\bone = \bone$ by construction, it is usual to expect that~$E$ is not precisely
the zero matrix.  
It is shown in Appendix~\ref{NE-details} that, if the quantity
\[
     \lmax\{E\TT E\} + \Tr\{W^2 + W W\TT\} \Eq \lmax\{(WX)\TT (I-H)WX\} + \Tr\{W^2 + W W\TT\}
\]
is not large, 
then neither $\rho$ nor~$\b$ can be accurately estimated.
}

\msk\nin\adbb{{\bf Conclusion~6:}  In the network effects model,
if the columns of~$WX$ (almost) belong to the space generated by the columns of~$X$, then estimating~$\r$ 
and the parameters~$\b$ in uniformly dense networks may both be problematic.
}

\section{Conclusion}\label{conclusion}
\setcounter{equation}{0}
\adbb{
Maximum likelihood estimation of the correlation parameter~$\r$ in the network disturbance and network effects
models suffers from substantial bias and variability when the underlying graph is dense, as has been
documented, for instance, in  \cite{mizruchi2008effect},  \cite{smith2009estimation}, \cite{farber2010topology},  
\cite{neuman2010structure} and \cite{la2018finite}.   In the context of the network disturbance model~\Ref{M-rho-def}, 
we have shown
that, in dense networks, there is a lower bound to the variability of {\it any\/} reasonable
estimator of~$\rho$, as a consequence of Cram\'er--Rao theory, and that increasing the size
of the network does not reduce the variability, if the density remains constant.  We also demonstrate that, 
under similar circumstances, the maximum likelihood estimator is significantly biased.  For certain ill conditioned
networks, these effects are particularly marked.  In Section~\ref{new-estimator}, we
propose a simple estimator of~$\rho$ that has little bias, though its variability in dense networks is
still substantial, because all estimators are subject to the limitations imposed by the Cram\'er--Rao 
lower bound.  We also suggest a permutation algorithm for assessing the variability of the estimator,
on the basis of data consisting of a single network.  The effectiveness of the new estimator is illustrated
by simulation 
in a few simple settings in Section~\ref{simul}.  There is a discussion of network effects 
models in Section~\ref{NE-models-sect}.
}
\adb{Here, in circumstances in which~$\r$ cannot be accurately
estimated, the estimates of the structural parameters~$\b$ are also affected.  An explicit quantity is given that 
indicates whether such problems are likely to occur.  It is expressed in terms 
of the design matrix~$X$ and the weight matrix~$W$.
}

\section*{Acknowledgement}
ADB thanks the mathematics departments of the University of Melbourne and Monash University, for
their kind hospitality while part of the work was undertaken. 
GR was supported in part by EP/T018445/1 and EP/R018472/1.
Both authors express their warm appreciation of the many helpful suggestions by the referees; these
have greatly improved the presentation.

\section*{Appendix}\label{Appendix}
\renewcommand{\thesubsection}{A.\arabic{subsection}}
\renewcommand{\theequation}{A.\arabic{equation}}
\setcounter{equation}{0}

\subsection{Rewriting the likelihood equation for~$\r$}\label{U-calc}
To minimize the quantity~$F(\r;y)$ of~\Ref{infer-disturbance-rho-hat}, it is usual to look for values of~$\r$ 
such that $\frac{d}{d\r} F(\r;y) =0$. First, as for~\Ref{urho-def},
\eq\label{D-det}
    \frac d{d\r}\,\log\det K_\r \Eq - \Tr\{Z_\r\}.
\en
Then the derivative of~$\shat_\r^2$ in~\Ref{NDM-MLE} is given by
\eq\label{sigma2-deriv}
   n\, \frac{d\shat^2_\r(y)}{d\r} 
          \Eq y\TT\Bigl\{-W\TT(I-H_\r)K_\r - K_\r\TT\Bigl(\frac{dH_\r}{d\r}\Bigr)K_\r - K_\r\TT(I-H_\r)W\Bigr\}y.
\en
To compute $\frac{d}{d\r}H_\r$, note that,  \adbb{for any family of invertible matrices $(Y_x,\,x \in \re)$ with
differentiable components, $0 = \frac d{dx} \{Y_x^{-1} Y_x\}$, implying that
\eq\label{Matrix-inverse-deriv}
     \frac d{dx} \{Y_x^{-1}\}  \Eq - Y_x^{-1} \frac d{dx} \{Y_x\} Y_x^{-1}.
\en 
Hence,} writing $Y_\r := ((K_\r X)\TT K_\r X)$, \adbb{and using $\frac{dK_\r}{d\r} = -W$,} we have
\eqs
      \frac{ dY_\r^{-1}}{d\r} 
          &=& Y_\r^{-1}\{X\TT W\TT K_\r X + X\TT K_\r\TT W X\} Y_\r^{-1}.
\ens
\adbb{Thus, noting for the last line that $W = Z_\r K_\r$ and that $K_\r X Y_\r^{-1}(K_\r X)\TT = H_\r$,} we deduce that
\eqa
    \frac {dH_\r}{d\r} &=& \frac d{d\r} \Bigl\{K_\r X\, Y_\r^{-1} (K_\r X)\TT \Bigr\} \non\\
             &=& -WX Y_\r^{-1}(K_\r X)\TT - K_\r X Y_\r^{-1}(WX)\TT \non\\
              &&\qquad\mbox{}  + K_\r X Y_\r^{-1}\{X\TT W\TT K_\r X + X\TT K_\r\TT W X\}Y_\r^{-1} (K_\r X)\TT \non\\
             &=& - (I-H_\r)Z_\r H_\r - H_\r Z_\r\TT(I-H_\r).  \label{DH-formula}
\ena
Note also that, because $H_\r$ is idempotent, which implies that
$(I-H_\r)H_\r = 0$, we have
\eq\label{mean-of-DH}
     \Tr((Z_\r H_\r)\TT(I-H_\r)) \Eq \Tr((I-H_\r)Z_\r H_\r) \Eq \Tr(Z_\r H_\r(I-H_\r)) \Eq 0 .
\en
Combining  \Ref{sigma2-deriv} and~\Ref{DH-formula}, it follows that
\eq\label{sigma-deriv-new}
   n\, \frac{d\shat^2_\r(y)}{d\r} \Eq - (K_\r y)\TT (I - H_\r)(Z_\r\TT + Z_\r)(I-H_\r)  K_\r y,
\en
and thus, writing $U(\r;y) \Def \frac{d}{d\r}F(\r;y)$ and using~\Ref{D-det},
\eqa
 U(\r;y) &=& \frac1{\shat^2_\r(y)}\,\frac{d\shat^2_\r(y)}{d\r}  - \frac2n\, \frac d{d\r}\,\log\det K_\r 
               \label{D-F-a}\\
    &=& -\frac1{n\shat^2_\r(y)}\, (K_\r y)\TT (I - H_\r)(Z_\r\TT + Z_\r)(I-H_\r) K_\r y  
           + \frac2n \sum_{l \ge 1}\r^l\Tr\{W^{l+1}\}, \phantom{XX} \label{D-F-1}
\ena
as given in~\Ref{D-F}.

\subsection{The magnitude of $nU'(\r;y)$}\label{U-deriv}

From \Ref{D-F-a} and~\Ref{D-F-1}, and since, for any $n\times n$ matrix~$M$ and $n$-vector $y$, $y\TT M y = y\TT M\TT y$,
it follows that
\eqa
     U'(\r;y) &=&  \frac1{\shat^2_\r(y)}\,\frac{d^2\shat^2_\r(y)}{d\r^2} 
            - \Bigl\{\frac1{\shat^2_\r(y)}\,\frac{d\shat^2_\r(y)}{d\r}\Bigr\}^2 + \frac2n\sum_{l\ge1}l\r^{l-1}\Tr\{W^{l+1}\}
                         \non\\
        &=& -\frac2{n\shat^2_\r(y)} \frac d{d\r}\bigl\{y\TT K_{\r}\TT\, \{(I-H_{\r})Z_\r (I-H_{\r})\}\,K_\r y\bigr\} \non\\
         &&\qquad\mbox{}      - \Bigl\{ U(\r;y) - \frac2n \sum_{l \ge 1}\r^l\Tr\{W^{l+1}\} \Bigr\}^2
                + \frac2n\sum_{l\ge1}l\r^{l-1}\Tr\{W^{l+1}\},\phantom{XX}\label{nU-again}
\ena
invoking~\Ref{sigma-deriv-new} for the first term.
Using $\frac d{d\r}K_\r = -W \adbb{= -Z_\r K_\r}$ and 
 \adbb{$\frac d{d\r}Z_\r = W \frac d{d\r}\{K_\r^{-1}\} = Z_\r^2$, from~\Ref{Matrix-inverse-deriv},}
and again using $y\TT M y = y\TT M\TT y$,
we can express the derivative of the quadratic form in~\Ref{nU-again} in terms of the quantities 
 $H_\r$, $K_\r$, $Z_\r$ and~$\frac d{d\r}H_\r$ as
\eqa
     \lefteqn{ \frac d{d\r}\{y\TT K_{\r}\TT\, \{(I-H_{\r})Z_\r (I-H_{\r})\}\,K_\r y\} }  \non\\
     &=&  y\TT K_\r\TT \Bigl\{ -Z_\r\TT (I-H_{\r})(Z_\r+Z_\r\TT) (I-H_{\r}) 
                               -  \frac{dH_\r}{d\r}(Z_\r+Z_\r\TT) (I-H_{\r}) \non \\
     &&\qquad\mbox{} +  (I-H_{\r})Z_\r^2 (I-H_{\r}) \Bigr\}K_\r y. \label{U-dash-1}
\ena
The middle term in~\Ref{nU-again} is expressed in terms of the same quantities
\adb{$H_\r$, $K_\r$, $Z_\r$ and~$\frac d{d\r}H_\r$,} together with~$W$, using~\Ref{sigma2-deriv}.
The final term in~\Ref{nU-again}, much as for~\Ref{trace-Z_rho}, gives~$2n^{-1}\Tr\{Z_\r^2\}$, and, 
from~\Ref{DH-formula},
\eqa
    \lefteqn{ y\TT K_\r\TT \frac{dH_\r}{d\r}(Z_\r+Z_\r\TT)(I-H_{\r})K_\r y }\non\\
      && \Eq - y\TT K_\r\TT \{(I-H_\r)Z_\r H_\r + H_\r Z_\r\TT(I - H_\r)\}(Z_\r + Z_\r\TT)(I-H_{\r})K_\r y.
            \label{U-dash-2}
\ena
The resulting expressions for~$U'(\r;y)$ are not easy to understand.  

It is somewhat easier to consider, as a typical value, the expectation
$n\ex_{\r_0} U'(\r_0;\by)$ on the model~$M_{\r_0}$. Approximating $\shat^2_\r(\by)$ by~$\s^2$, the expectation of~$n$
times the first and third terms in~\Ref{nU-again} is approximated by
\[
   \ex_{\r_0}\Bigl\{ -\frac2{\s^2} \frac d{d\r}\{\by\TT K_{\r}\TT\, \{(I-H_{\r})Z_\r (I-H_{\r})\}\,K_\r \by\} 
                          \adb{\Bigr|_{\r=\r_0}}\Bigr\}
              + 2 \Tr\{Z_{\r_0}^2\}.
\]
\adbb{
This expression, using~\Ref{quad-form-mean} together with \Ref{U-dash-1} and~\Ref{U-dash-2} \adb{and properties of the trace, 
and also observing that $(I - H_\r)K_\r X = 0$,} gives
\eqa
    \lefteqn{2\Tr\{Z_{\r_0}\TT (I-H_{{\r_0}})(Z_{\r_0}+Z_{\r_0}\TT) (I-H_{{\r_0}}) } \non\\  
     &&\qquad\mbox{} - \{(I-H_{\r_0})Z_{\r_0} H_{\r_0} + H_{\r_0} Z_{\r_0}\TT(I - H_{\r_0})\}
                              (Z_{\r_0} + Z_{\r_0}\TT)(I-H_{{\r_0}}) \phantom{XXXX}\non\\
     &&\qquad\qquad\mbox{} - (I-H_{{\r_0}})Z_{\r_0}^2 (I-H_{{\r_0}}) + Z_{\r_0}^2\} \non\\
     &&\Eq \Tr\{(Z_{\r_0} + Z_{\r_0}\TT) (I-H_{{\r_0}})(Z_{\r_0}+Z_{\r_0}\TT) (I-H_{{\r_0}})\} \non\\
     &&\qquad\mbox{}  + 2\Tr\{ - Z_{\r_0} H_{\r_0}(Z_{\r_0} + Z_{\r_0}\TT)(I - H_{\r_0})  +  H_{\r_0} Z_{\r_0}^2\} 
                 \ =:\ \D_0,
               \label{U-dash-expec-1}
\ena
with the first trace equal to~$\half\t_0^2$, \adb{ as given in~\Ref{tau0-def}.}
}
The middle term in~\Ref{nU-again}, multiplied by~$n$, gives
\[
    n^{-1}\var_{\r_0}(nU(\r;\by)) + n^{-1}(n\ex_{\r_0}U(\r;\by) - 2\Tr\{Z_0\})^2 \ \approx\
        n^{-1}(\t_0^2  + 4\Tr\{(I-H_{\r_0})Z_{\r_0}\}^2),     
\]
using~\Ref{mu0-def} and~\Ref{tau0-def}.  In the \adb{ circumstances discussed in Section~\ref{Estn-with-structure},} 
$\D_0$ and~$\t_0^2$ are seen to be
typically comparable to $\Tr\{Z_{\r_0}Z_{\r_0}\TT\}$, so that this last contribution is of relative
order~$O(n^{-1})$ as compared to~$\D_0$.  Hence, in such circumstances, the expectation
$n\ex_{\r_0} U'(\r_0;\by)$ on the model~$M_{\r_0}$ is approximately given by~$\D_0$.

\subsection{Details of the MLE when structure is present}\label{4.3-extras}
Suppose that~$X$ consists solely of the $n$-vector $\bone$, having all elements equal to~$1$, and that 
there are no isolated vertices. 
Suppose also that~$W$ is chosen to have $\sjn W_{ij} = 1$ for all~$j$, so that $W\bone = \bone$. 
Then $H_{\r} = H_* := n^{-1}\bone\bone\TT$ for all~$\r$,
because $K_\r X = (I-\r W)\bone = (1-\r)\bone$, and hence $WH_\r = H_*$ and $Z_\r H_{\r} =  (1-\r)^{-1}H_*$.
\adb{
The matrix product in~\Ref{tau0-def} can be multiplied out to give
\eqa
    \t_0^2 &=& 2 \Tr\{(Z_{\r_0} + Z_{\r_0}\TT) (I - H_{\r_0}) (Z_{\r_0} + Z_{\r_0}\TT) (I - H_{\r_0})\} \non\\
           &=& 2 \Tr\{ 2(Z_{\r_0}^2 + Z_{\r_0}Z_{\r_0}\TT) -  2H_{\r_0} (Z_{\r_0} + Z_{\r_0}\TT)^2
                          + [H_{\r_0} (Z_{\r_0} + Z_{\r_0}\TT)]^2 \} \label{tau-0-expanded}\\
           &=& 4\Tr\{ (Z_{\r_0}^2 + Z_{\r_0}Z_{\r_0}\TT)(I - H_{\r_0})\}
                - 2\Tr\{H_{\r_0} (Z_{\r_0} + Z_{\r_0}\TT)^2 - [H_{\r_0} (Z_{\r_0} + Z_{\r_0}\TT)]^2 \}.
                     \non
\ena
With the particular choices of $X$ and~$W$, and repeatedly using the  properties of the trace, this gives
\eqs
     \t_0^2 
            &=& 4\Tr\{ (Z_{\r_0}^2 + Z_{\r_0}Z_{\r_0}\TT)(I - H_*) \} + 
                      \red{2}\Bigl\{\frac1{(1-\r_0)^2} - \frac1n |Z_{\r_0}\TT\bone|^2 \Bigr\},
\ens
where, for instance, the final term is just $-\Tr\{H_*Z_{\r_0} Z_{\r_0}\TT\}$, whereas, for example, 
$$
      \Tr\{H_*Z_{\r_0}\TT Z_{\r_0}\} \Eq \frac1{1-\r_0}\Tr\{H_* Z_{\r_0}\} \Eq \frac1{1-\r_0}\Tr\{Z_{\r_0}H_*\}
                    \Eq \frac1{(1-\r_0)^2}\,.
$$
In much the same way, and in particular because $Z_\r H_* =  (1-\r)^{-1}H_*$ and $H_*(I - H_*) = 0$, we obtain
\[
 2\Tr\{ - Z_{\r_0} H_*(Z_{\r_0} + Z_{\r_0}\TT)(I - H_*)  +  H_* Z_{\r_0}^2\} \Eq 2\Tr\{Z_{\r_0}^2 H_*\}
               \Eq \frac2{(1-\r_0)^2}\,,
\]
giving the formula for~$\D_0$ in~\Ref{with-structure}.}

\adb{Note that 
$$
         \frac1n |Z_{\r_0}\TT\bone|^2\ \ge\ \Bigl\{  \frac1n \bone\TT Z_{\r_0}\TT\bone \Bigr\}^2 \Eq \frac1{(1-\r_0)^2}\,,
$$
by \red{Cauchy--Schwarz,} giving the upper bound for~$\t_0^2$ in~\Ref{with-structure}.}
Note also that, again using $Z_\r H_{\r} =  (1-\r)^{-1}H_*$, and because $\Tr\{H_*\} = 1$ \adb{and $H_*^2 = H_*$}, 
it follows that
\eqa
     \quarter \t_0^2 &=& \Tr\{Z_{\r_0}Z_{\r_0}\TT + Z_{\r_0}^2\} - \Tr\{H_* Z_{\r_0}(H_* Z_{\r_0})\TT\} 
                        - \Tr\{Z_{\r_0}^2 H_*\}  \non\\
                     &\le& \Tr\{Z_{\r_0}Z_{\r_0}\TT + Z_{\r_0}^2\} - \frac1{(1-\r_0)^2}\,, \label{quarter-tau-squared}
\ena
so that $\t_0^2$ is indeed smaller than $4\Tr\{Z_{\r_0}Z_{\r_0}\TT + Z_{\r_0}^2\}$.

\adb{
Suppose now that there is more structure in the regression, so that~$X$ consists of columns in addition to the
first column~$\bone$.  If we still assume} that $W\bone = \bone$, the basic message remains the same.
Because the matrix~$H_\r$ is idempotent, its spectral decomposition can be written as
\eq\label{H-spectral}
     H_\r \Eq H_* + \sum_{l=2}^m x_\r\ul (x_\r\ul)\TT,
\en
where $H_* := n^{-1}\bone\bone\TT$ is as before, and $x_\r\ut,\ldots,x_\r\um$ are {\it orthonormal\/} vectors 
orthogonal to~$\bone$ that, with~$\bone$, span the space generated by the columns of $K_\r X$.
If, for instance, it is assumed that each of the vectors~$x_\r\ul$,
$2\le l\le m$, has components bounded in modulus by~$cn^{-1/2}$, for some fixed~$c \ge 1$ --- 
as is the case for $x\ui = n^{-1/2}\bone$ --- then the traces of the matrices $Z_{\r_0}^2 H_{\r_0}$, 
$Z_{\r_0}\TT Z_{\r_0} H_{\r_0}$, $(Z_{\r_0} H_{\r_0})^2$ and $Z_{\r_0}\TT H_{\r_0} Z_{\r_0} H_{\r_0}$
\adb{can all be shown to be bounded by  $m^2 c^4/(1-|\r|)^2$, which is not large if $m$ and~$c$ are moderate
and~$|\r|$ is not too close to~$1$.   The argument runs as follows.
}

In view of~\Ref{H-rho-general}, and writing $x\ui := n^{-1/2}\bone$,  the trace $\Tr\{MH_\r\}$,
for any $n\times n$ matrix~$M$, can be written in the form
\eq\label{trace-MH}
     \Tr\{MH_\r\} \Eq \sum_{l=1}^m \Tr\{Mx_\r\ul(x_\r\ul)\TT\}  \Eq \sum_{l=1}^m (x_\r\ul)\TT M x_\r\ul.
\en
Thus, using~\Ref{eps-from-nu-expansion}, it follows that, for any $1 \le l,l' \le m$, 
\adb{
\eqa
   |Z_\r x_\r\ul| &=& \Bigl| \sum_{r\ge0}  W^{r+1} \r^r x_\r\ul \Bigr| \non\\
           &\le& \sum_{r\ge0}|\r|^r  cn^{-1/2}|W^{r+1} \bone| \Eq \frac{c}{n^{1/2}(1-|\r|)}|\bone|
                   \Eq  \frac{c}{1 - |\r|}\,; \label{x-trace-0}
\ena
}
that
\eqa
    |(x_\r\ul)\TT Z_\r x_\r\uld| &=& \Bigl| \sum_{r\ge0} (x_\r\ul)\TT W^{r+1} \r^r x_\r\uld \Bigr| \non\\
      &\le& \adbb{\Bigl| \sum_{r\ge0}|\r|^r (cn^{-1/2}\bone\TT)  W^{r+1} (cn^{-1/2}\bone) \Bigr|}
    \Eq c^2n^{-1} \adbb{\bone\TT\bone} \sum_{r\ge0}  |\r|^r  \non\\
      &=& \frac{c^2}{1 - |\r|}\,; \label{x-trace-1}
\ena
and that
\adbb{
\eqa
    |(x_\r\ul)\TT Z_\r^2 x_\r\uld| &=& |(x_\r\ul)\TT W^2 (1 - \r W)^{-2} x_\r\uld| 
         \Eq \Bigl|  \sum_{r\ge0} (x_\r\ul)\TT (r+1)W^{r+2} \r^r x_\r\uld \Bigr|  \non\\ 
        &\le& c^2 n^{-1} \bone\TT\bone \sum_{r\ge0} (r+1)|\r|^r \Eq \frac{c^2}{(1 - |\r|)^2}\,. \label{x-trace-2}
\ena
}
\adb{
Then, for instance, from~\Ref{trace-MH} and~\Ref{x-trace-2}, we have
\eqs
    \Tr\{Z_{\r_0}^2 H_{\r_0}\} &=& \sum_{l=1}^m (x_\r\ul)\TT Z_{\r_0}^2 x_\r\ul 
             \Le  \frac{mc^2}{(1 - |\r|)^2}\,;
\ens
from~\Ref{trace-MH} and~\Ref{x-trace-0}, we have
\eqs
    \Tr\{Z_{\r_0}\TT Z_{\r_0} H_{\r_0}\} &=& \sum_{l=1}^m (x_\r\ul)\TT Z_{\r_0}\TT Z_{\r_0} x_\r\ul 
             \Eq \sum_{l=1}^m |Z_{\r_0} x_\r\ul|^2 \Le \frac{mc^2}{(1 - |\r|)^2}\,;
\ens
and, from~\Ref{trace-MH} and~\Ref{x-trace-1}, we have
\eqs
  \Tr\{Z_{\r_0}\TT H_{\r_0} Z_{\r_0} H_{\r_0}\} &=& \sum_{l=1}^m (x_\r\ul)\TT Z_{\r_0}\TT H_{\r_0} Z_{\r_0} x_\r\ul 
      \Eq \sum_{l=1}^m \Tr\{Z_{\r_0} x_\r\ul (x_\r\ul)\TT Z_{\r_0}\TT H_{\r_0} \} \\
      &=& \sum_{l=1}^m \sum_{l'=1}^m (x_\r\uld)\TT Z_{\r_0} x_\r\ul (x_\r\ul)\TT Z_{\r_0}\TT  x_\r\uld 
      \Le \frac{m^2c^4}{(1 - |\r|)^2}\,.
\ens
By similar calculations, each of the traces listed above is bounded in modulus by the quantity $m^2 c^4/(1-|\r|)^2$.
}

\subsection{Estimation of~$\b$}\label{off-the-model}
\adb{
From \Ref{NDM-MLE} and~\Ref{NDM-MLE-2},
the distribution of~$\hat\b_\r$ on the model~$M_\r$ has covariance matrix $\s^2 (X\TT\VV_\r X)^{-1}$,
where $\VV_\r := K_\r \TT K_\r$.  The variance of a linear combination~$a\TT \hat\b_\r$, 
for a given $m$-dimensional unit vector~$a$, is thus
\eq\label{off-model-1}
     \s^2 a\TT (X\TT\VV_\r X)^{-1} a \Le \s^2 \lmax\{(X\TT\VV_\r X)^{-1}\}
          \Eq \s^2 / \lmin(X\TT\VV_\r X),
\en
where $\lmin(M) \le \lmax(M)$ denote the smallest and largest eigenvalues of~$M$, both positive if~$M$
is positive definite symmetric.   
Letting $\tu_X$ denote the unit $n$-vector $Xu / |Xu|$, and using Rayleigh--Ritz twice, we have
\eqa
   \lmin(X\TT\VV_\r X) &=& \min_{u\in \re^m \colon |u|=1} u\TT X\TT\VV_\r X u 
                       \Eq \min_{u\in \re^m \colon |u|=1} |Xu|^2 \tu_X\TT \VV_\r  \tu_X  \non\\
              &\ge& \lmin(\VV_\r) \min_{u\in \re^m \colon |u|=1} |Xu|^2 \Eq \lmin(X\TT X) \lmin(\VV_\r). \label{off-model-2}
\ena
Hence, for fixed~$\s^2$, provided that $\lmin(\VV_\r)$ is not close to zero, \Ref{off-model-1}
and~\Ref{off-model-2} imply that accurate estimation of
all linear combinations of~$\b$ depends on having a large value of $\lmin(X\TT X)$, which is the usual
condition for ordinary least squares.  In particular, if the network size~$n$ becomes large, but
$\lmin(\VV_\r)$ remains uniformly bounded away from zero, the estimate of~$\b$ is consistent if
$\lmin(X\TT X) \to \infty$. \red{Example~1 in Section~\ref{Ill-conditioned} shows that the latter condition
alone is not enough to guarantee consistent estimation of~$\b$.}
}

\adb{
If the true value of $\r$ is~$\r_0$, the distribution of~$\hat\b_\r(\by)$, as determined by substituting
$\by = X\b + K_{\r_0}^{-1}\bnu$ for~$y$ into~\Ref{NDM-MLE}, is multivariate normal, with mean~$\b$ and covariance
matrix $\s^2\Sigma(\r,\r_0)$, where
\eqa
   \Sigma(\r,\r_0) &:=& 
         (X\TT \VV_\r X)^{-1}(K_\r X)\TT K_\r K_{\r_0}^{-1} (K_{\r_0}^{-1})\TT K_\r\TT  K_\r X (X\TT \VV_\r X)^{-1} \non\\
    &=& (X\TT \VV_\r X)^{-1} X\TT \VV_\r  (\VV_{\r_0})^{-1} \VV_\r X (X\TT \VV_\r X)^{-1}. \label{off-model-3}
\ena
Writing $\tX_\r := \VV_{\r}^{1/2}X(X\TT \VV_\r X)^{-1}$, where $M^{1/2}$ denotes the non-negative square root of 
a non-negative definite
matrix~$M$, and setting $\VV(\r,\r_0) := \VV_{\r}^{1/2} (\VV_{\r_0})^{-1} \VV_{\r}^{1/2}$, \Ref{off-model-3} gives
\eq\label{off-model-4}
      \Sigma(\r,\r_0)  \Eq   \tX_\r\TT \VV(\r,\r_0) \tX_\r .
\en
Using Rayleigh--Ritz, for any unit $m$-vector $a$, we thus have
\eqs
    a\TT  \Sigma(\r,\r_0) a \Eq (\tX_\r a)\TT \VV(\r,\r_0) \tX_\r a
                 &\le& |\tX_\r a|^2 \lmax(\VV(\r,\r_0)) \\
                 &=& \{a\TT (X\TT \VV_\r X)^{-1} a\} \lmax(\VV(\r,\r_0)).
\ens
Thus, if the true model is~$M_{\r_0}$, estimating a linear combination~$a\TT\b$ using $\hat\b_{\r}$ instead 
of~$\hat\b_{\r_0}$ gives an estimator  $a\TT\hat\b_{\r}$ whose variance is at most $\lmax(\VV(\r,\r_0))$
times the value, \red{given in~\Ref{off-model-1},} that it would have had, were~$M_\r$ the true model.  
Hence, if estimation on the model~$M_\r$,
using~$\hat\b_\r$, is accurate for any~$\r$, then estimating using~$\hat\b_\r$ on the model~$M_{\r_0}$ is also
accurate, provided that $\lmax(\VV(\r,\r_0))$ is not large.
} 

\red{Note that, if $n(W) := \sqrt{\lmax(W\TT W)}$ denotes the operator norm of the matrix~$W$, then, 
by Rayleigh--Ritz, \adbp{because $|Wy| \le n(W)|y|$ for all~$y$,} the largest eigenvalue 
of $S(\r,0) = I - \r(W + W\TT) + \r^2 W\TT W$ is at least $(|\r|n(W) - 1)^2$, and this
is large, for moderate non-zero values of~$\r$, if~$n(W)$ is large.  Hence, if~$\lmax(W\TT W)$ is large, this
indicates that there may possibly be problems with estimating~$\b$.}

\subsection{The scale of precision of estimation using $T^C_\r$}\label{scale-quadratic-estimator}
As discussed in Section~\ref{Quadratic-form-estimation},
a plausible measure of the scale of precision of estimation using $T^C_\r$ is given by
$\ssch^C_{\r_0} := \htau^C_{\r_0}/|\hD^C_{\r_0}|$, to be estimated by $\htau^C_{\hrho_C}/|\hD^C_{\hrho_C}| $, where
\[
    (\htau^C_{\r_0})^2 \Def \s^{-4}\var_{\r_0}(U^C_{\r_0}(\by)) \quad\mbox{and}\quad 
                \hD^C_{\r_0} \Def \adb{-} \s^{-2}\ex_{\r_0}\Bigl\{\frac{dU^C_{\r}(\by)}{d\r}\Bigr|_{\r=\r_0}\Bigr\},
\]
\adb{and $U^C$ is as defined in~\Ref{U^C-def}.}

\adbb{To find~$(\htau^C_{\r_0})^2$, note that, from \Ref{A-QF} and~\Ref{NDM-MLE}, since 
$\by = X\b + K_{\r_0}^{-1}\bnu$ under~$M_{\r_0}$, we have
\eqs
    U^C_{\r_0}(\by) &=& \bnu\TT (Q_C(\r_0,\r_0)  + n^{-1}(I - H_{\r_0})\Tr\{H_{\r_0}C\})\bnu \Eq \bnu\TT \hQ(C,\r_0)\bnu,
\ens
\adb{where $Q_C(\r,\r_0)$ is defined in~\Ref{QA-def},} and
\eq\label{QC-def}
     \hQ(C,\r_0)  \Def (I-H_{\r_0}) C (I-H_{\r_0}) + n^{-1}(I - H_{\r_0})\Tr\{H_{\r_0}C\}.
\en
Hence, replacing $C$ by $\half(C + C\TT)$, \adb{which does not change~$T^C$,} and using~\Ref{quad-form-var}, 
it follows that
\eqa\label{tau-hat-def}
    (\htau^C_{\r_0})^2 &=& \s^{-4}\var_{\r_0}( U^C_{\r_0}(\by)) \Eq \half \Tr\bigl\{\{\hQ(C + C\TT,\r_0)\}^2\bigr\}.
\ena
In view of the factor~$n^{-1}$ in the second element of~$\hQ(C,\r_0)$, the dominant contribution is that from 
\[
  \half \Tr\bigl\{\{Q_C(\r_0,\r_0)\}^2\bigr\} \Eq \Tr\{(I-H_{\r_0}) C (I-H_{\r_0})(C + C\TT)\},
\]
once more of the form $\Tr\{Z^2 + ZZ\TT\}$, now with $Z := (I-H_{\r_0})C(I-H_{\r_0})$.
}

\adbb{
To find $\hD^C_{\r_0}$,} observe that, since $\frac d{d\r}K_\r = -W = -Z_{\r}K_{\r}$ and $\frac d{d\r}H_\r$ is 
as in~\Ref{DH-formula}, we have
\eqa
    \frac{dT^C_\r(\by)}{d\r} &=& \frac{d}{d\r}\{(K_\r\by)\TT (I-H_\r) C (I - H_\r) K_\r \by\} \non\\
        &=& - (K_{\r}\by)\TT \Bigl( Z_\r\TT(I - H_\r) + \frac{dH_\r}{d\r} \Bigr)(C + C\TT) (I-H_\r) K_\r \by
                     \phantom{XXXXXXXXXX} \non\\
        &=& - (\adbb{X\b +} K_{\r_0}^{-1}\bnu)\TT K_\r\TT  \bigl\{Z_\r\TT (I - H_\r)  - (I-H_\r)Z_\r H_\r \non\\
        &&\qquad\qquad\qquad\qquad \mbox{} - H_\r Z_\r\TT(I-H_\r) \bigr\} (C + C\TT) (I-H_\r) K_\r K_{\r_0}^{-1}\bnu,
                   \label{T-deriv}  
\ena
\adbb{where~$X\b$ is not present at the second appearance of~$\bnu$ because $(I-H_\r) K_\r X = 0$.}
 From \Ref{T-deriv} \adbb{and~\Ref{quad-form-mean}, and because $H_\r(I-H_\r) = 0$,} it thus follows that
\eq\label{Delta-hat-1-def}
    -\ex_{\r_0}\Bigl\{\frac{dT^C_\r(\by)}{d\r}\Big|_{\r = \r_0}\Bigr\}
             \Eq \s^2\Tr\bigl\{(I-H_{\r_0})\adbb{[Z_{\r_0}\TT(I - H_{\r_0}) - Z_{\r_0}H_{\r_0}]} (C + C\TT) \bigr\}
             \ =:\ \s^2\hD^C_{\r_0,1} ,
\en
\adbb{say, with~$X\b$ disappearing because $\ex \bnu = 0$.  Then, from~\Ref{DH-formula},
\eq\label{Delta-hat-2-def}
     -\frac{d}{d\r} \Tr\{H_\r C\} \Eq -\Tr\Bigl\{C\,\frac{dH_\r}{d\r}\Bigr\} \Eq
             \Tr\bigl\{C\{(I-H_\r)Z_\r H_\r + H_\r Z_\r\TT(I-H_\r)\}\bigr\} \ =:\ \hD^C_{\r_0,2},
\en
say, and, from~\Ref{sigma-deriv-new}, and on~$M_{\r_0}$,
\[
  -\frac{d\shat^2_{\r_0}(\by)}{d\r} \Eq  n^{-1} \bnu\TT (I - H_{\r_0})(Z_{\r_0}\TT + Z_{\r_0})(I-H_{\r_0}) \bnu,
\]
so that 
\eq\label{Delta-hat-3-def}
    -\s^{-2}\ex_{\r_0}\Bigl\{ \frac{d\shat^2_{\r_0}(\by)}{d\r} \Bigr\} \Tr\{H_{\r_0}C\}
          \Eq 2n^{-1} \Tr\{(I - H_{\r_0})Z_{\r_0}(I-H_{\r_0})\} \Tr\{H_{\r_0}C\} \ =:\ \hD^C_{\r_0,3},
\en
say.  Then 
\eq\label{Delta-hat-def}
     \hD^C_{\r_0} \Def \sum_{l=1}^3 \hD^C_{\r_0,l},
\en
where $\hD^C_{\r_0,l}$, $1\le l\le 3$, are given by \Ref{Delta-hat-1-def}, \Ref{Delta-hat-2-def} and~\Ref{Delta-hat-3-def}.
}

\subsection{Details of network effects estimation}\label{NE-details}
\adbb{
In the network effects model, except when only the overall mean is to be estimated,
and $W\bone = \bone$ by construction, it is usual to expect that
\eq\label{WX-decomposition}
     WX \Eq HWX + (I-H)WX \Eq X\G + E,
\en
for a non-zero $n \times m$ matrix~$E$.  Here, $H := X(X\TT X)^{-1}X\TT$ is the projection onto the space spanned by 
the columns of~$X$; $HWX$ is the corresponding projection of the columns of~$WX$, and can thus be written in the 
form~$X\G$; and $E := (I-H)WX$ has columns orthogonal to the columns of~$X$.  Then, as in the calculation leading
to~\Ref{KXb-rephrase}, we obtain
\eqs
   K_\r^{-1}X\b &=&  X(I - \r\G)^{-1}\b + \r K_\r^{-1} E \b \\
                &=&  X\parg + \r K_\r^{-1} E (I - \r\G) \parg \ =:\ X_\r \parg,
\ens
say, parametrized in terms of~$\parg$.  The log-likelihood, in terms of~$\parg$, becomes
\[
     \ell(\parg,\s^2,\r;y) \Eq 
       -\tfrac n2 \{\log(\s^2) + \log (2\pi)\} - \frac1{2\s^2} \{K_\r(y-X_\r\parg)\}\TT  \{K_\r (y-X_\r\parg)\} 
               + \log\det K_\r,
\]
and differentiating with respect to~$\r$ gives the analogue of the score function in~\Ref{urho-def-2} to be
\eqa
   \ts\urho(\g,\s^2,\r;y) 
      &=& \frac1{\s^2} \bigl(\{K_\r(y-X_\r\parg)\}\TT  Z_\r \{K_\r(y-X_\r\parg)\}\bigr)  \non\\
       &&\qquad\mbox{}     + \frac1{\s^2} \{K_\r(y-X_\r\parg)\}\TT G_\r \parg  - \Tr\{Z_\r\}, \label{urho-def-3}
\ena
where
\eqa
     G_\r &:=& K_\r \frac d{d\r}\{\r K_\r^{-1} E (I - \r\G)\} \non\\
          &=& E (I - \r\G) - \r K_\r^{-1} WE (I - \r\G) - \r E \G.
\ena
The variance of $\ts\urho(\g_0,\s_0^2,\r_0;y)$ on the model with parameters $(\g_0,\s_0^2,\r_0)$ is now
\eq\label{var-score-2}
     V_2 \Def \Tr\{Z_{\r_0}^2 + Z_{\r_0} Z_{\r_0}\TT\} + \frac1{\s_0^2}|G_{\r_0}\parg|^2,
\en
differing from~\Ref{var-score-1} in the second term.  The corresponding Cram\'er--Rao lower bound for 
the variance of any unbiased estimator of~$\r$ is thus at least as big as~$1/V_2$.
}

\adbb{   Just as for the network
disturbance model, the first term in~\Ref{var-score-2} may not be large in dense networks. 
To get an idea of the magnitude of the second term, 
consider the case where $\r_0=0$, so that $G_0 = E$ and $|G_{0}\parg|^2 = \parg\TT E\TT E \parg$.
Then
\[
    \frac1{\s_0^2}|G_{0}\parg|^2 \Le \frac{|\parg|^2}{\s_0^2} \lmax(E\TT E),   
\]
where $\lmax(E\TT E)$ denotes the largest eigenvalue of the positive definite symmetric $m\times m$ matrix~$E\TT E$.
For instance, if the elements of~$E$, that represent to what extent the columns of~$WX$ do not belong
to the linear span of the columns of~$X$, are uniformly of magnitude at most $c/\sqrt n$, 
then the elements of $E\TT E$ are at most of magnitude~$c^2$, and thus~$\lmax(E\TT E)$ is not
large if~$c$ is not large.  In such circumstances, any reasonable
estimator of~$\r$ has to have substantial variability.  In practice,  the model decomposition~\Ref{WX-decomposition}
can be applied to yield the value of
\[
      \lmax\{E\TT E\} \Eq \lmax\{(WX)\TT (I-H)WX\},
\]
and, if neither this nor $\Tr\{W^2 + W W\TT\}$ is large, then neither $\r$ nor~$\b$ can be accurately estimated.
}

\end{document}